\documentclass[11pt, reqno]{amsart}

\usepackage[OT2, T1]{fontenc}

\usepackage[usenames,dvipsnames]{color}

\usepackage{url}
\usepackage{amsmath}
\usepackage{array}
\usepackage{graphicx}
\usepackage{amssymb}
\usepackage{amsthm}
\usepackage{colonequals}	% for nice :=
\usepackage{color}
\usepackage{mathabx}		% for \widec
\usepackage{mathrsfs}
\usepackage{amscd}
\usepackage[all,cmtip]{xy}	% for commutative diagrams
\usepackage{sseq}			% for spectral sequences
\usepackage{verbatim}
\usepackage{parskip}
\usepackage[in]{fullpage}
\usepackage{mathrsfs} 			% for \mathscr (script letters)
\usepackage{bm} 				% for bold Greek letters
\usepackage{enumitem}
\usepackage[alphabetic,lite]{amsrefs} 	% for bibliography
\usepackage{stmaryrd}	% this is for \llbracket and \rrbracket
\usepackage{hyperref}
\usepackage{cleveref}

\DeclareSymbolFont{cyrletters}{OT2}{wncyr}{m}{n}
\DeclareMathSymbol{\Sha}{\mathalpha}{cyrletters}{"58}

\usepackage{color}

\newcommand{\gA}{\alpha}

\newcommand{\bF}{\mathbb{F}}
\newcommand{\bG}{\mathbb{G}}

\newcommand{\bQ}{\mathbb{Q}}

\newcommand{\bZ}{\mathbb{Z}}

\newcommand{\cA}{\mathcal{A}}
\newcommand{\cB}{\mathcal{B}}

\newcommand{\cE}{\mathcal{E}}

\newcommand{\cG}{\mathcal{G}}

\newcommand{\cO}{\mathcal{O}}

\newcommand{\cT}{\mathcal{T}}

\newcommand{\cX}{\mathcal{X}}
\newcommand{\cY}{\mathcal{Y}}
\newcommand{\cZ}{\mathcal{Z}}

\newcommand{\fa}{\mathfrak{a}}

\newcommand{\fm}{\mathfrak{m}}
\newcommand{\fo}{\mathfrak{o}}
\newcommand{\fp}{\mathfrak{p}}

\newcommand{\ra}{\rightarrow}

\newcommand{\xra}{\xrightarrow}

\newcommand{\hra}{\hookrightarrow}

\newcommand{\I}{^{\infty}}
\newcommand{\wt}{\widetilde}
\newcommand{\wh}{\widehat}

\newcommand{\pr}{^{\prime}}
\newcommand{\ce}{\colonequals}

\renewcommand{\b}{\textbf}
\newcommand{\surjects}{\twoheadrightarrow}

\newcommand{\tensor}{\otimes} 		% binary tensor product
\newcommand{\isomto}{\overset{\sim}{\longrightarrow}}

\newcommand{\nr}{{\mathrm{nr}}}		% unramified cohomology or max unramified extension
\newcommand{\cris}{{\mathrm{cris}}}		% crystalline (mainly used in subscripts)
\newcommand{\fppf}{\mathrm{fppf}}		% for fppf cohomology (mainly used in subscripts)
\newcommand{\et}{\mathrm{\acute{e}t}}	% for etale cohomology (mainly used in subscripts)
\newcommand{\Et}{\mathrm{\acute{E}t}}	% for the big etale site (mainly used in subscripts)
\renewcommand{\implies}{\Rightarrow}
\providecommand{\p}[1]{\left(#1\right)}

\providecommand{\f}[2]{\frac{#1}{#2}}
\DeclareMathOperator{\Ker}{Ker}			% Kernel
\DeclareMathOperator{\im}{Im}			% Imaginary part
\DeclareMathOperator{\Spec}{Spec}		% Spectrum of a ring
\DeclareMathOperator{\Hom}{Hom}			% Set of arrows between two object
\DeclareMathOperator{\Char}{char}		% Characteristic of a field
\DeclareMathOperator{\Frac}{Frac}		% Field of fractions
\DeclareMathOperator{\Ext}{Ext}			% Derived functors of Hom
\DeclareMathOperator{\Gal}{Gal}	% Galois group
\DeclareMathOperator{\End}{End}		% The algebra of endomorphisms
\DeclareMathOperator{\rk}{rk}		% rank
\DeclareMathOperator{\cork}{cork}		% corank
\DeclareMathOperator{\Sel}{Sel}		% Selmer group
\DeclareMathOperator{\Lie}{Lie}		% Lie algebra
\DeclareMathOperator{\Cl}{Cl}		% Class group
\newcommand{\ba}{\begin{aligned}}
\newcommand{\ea}{\end{aligned}}
\newcommand{\be}{\begin{equation}}
\newcommand{\ee}{\end{equation}}
\newcommand{\pf}{\begin{proof}}
\newcommand{\bpf}{\begin{proof}}
\newcommand{\epf}{\end{proof}}
\newcommand{\bthm}{\begin{thm}}
\newcommand{\ethm}{\end{thm}}
\newcommand{\bprop}{\begin{prop}}
\newcommand{\eprop}{\end{prop}}
\newcommand{\bcor}{\begin{cor}}
\newcommand{\ecor}{\end{cor}}
\newcommand{\brem}{\begin{rem}}
\newcommand{\erem}{\end{rem}}
\newcommand{\brems}{\begin{rems} \hfill \begin{enumerate}[label=\b{\thesubsection.},ref=\thesubsection]}
\newcommand{\remi}{\addtocounter{subsection}{1} \item}
\newcommand{\erems}{\end{enumerate} \end{rems}}
\newcommand{\blem}{\begin{lemma}}
\newcommand{\elem}{\end{lemma}}
\newcommand{\bconj}{\begin{conj}}
\newcommand{\econj}{\end{conj}}
\newcommand{\benum}{\begin{enumerate}[label={(\alph*)}]}
\newcommand{\eenum}{\end{enumerate}}
\newcommand{\bc}{\begin{comment}}
\newcommand{\ec}{\end{comment}}
\newcommand{\beg}{\begin{eg}}
\newcommand{\eeg}{\end{eg}}
\newcommand{\lab}{\label}
\newcommand{\tst}{\textstyle}

\theoremstyle{plain}
\newtheorem{thm}[subsection]{Theorem}
\Crefname{thm}{Theorem}{Theorems}

\Crefname{rethm}{Theorem}{Theorem}

\Crefname{req}{Question}{Question}
\newtheorem{prop}[subsection]{Proposition}
\Crefname{prop}{Proposition}{Propositions}

\Crefname{q}{Question}{Questions}

\Crefname{Problem}{Problem}{Problems}
\newtheorem{conj}[subsection]{Conjecture}
\Crefname{conj}{Conjecture}{Conjectures}
\newtheorem{cor}[subsection]{Corollary}
\Crefname{cor}{Corollary}{Corollaries}

\newtheorem{lemma}[subsection]{Lemma}

\theoremstyle{definition}
\newtheorem{eg}[subsection]{Example}
\Crefname{eg}{Example}{Examples}

\theoremstyle{definition}
\newtheorem{rem}[subsection]{Remark}
\Crefname{rem}{Remark}{Remarks}
\newtheorem*{rems}{Remarks}
\Crefname{rems}{Remarks}{Remarks}

\theoremstyle{remark}

\Crefname{construction}{Construction}{Constructions}

\Crefname{claim}{Claim}{Claims}

\newtheoremstyle{subsection-tweak}
   {11pt}
   {3pt}%
   {}
   {}%
   {\bfseries}
   {}%
   {.5em}
   {\thmnumber{\@{#1}{}\@{#2}.}%
    \thmnote{~{\bfseries#3.}}}

\Crefname{innercustomconj}{Conjecture}{Conjecture}

\theoremstyle{subsection-tweak}
\newtheorem{pp}[subsection]{}
\newcommand{\bpp}{\begin{pp}}
\newcommand{\epp}{\end{pp}}

\makeatletter
\def\saveenum{\xdef\@savedenum{\the\c@subsection\relax}}
\def\resetenum{\global\c@enumi\@savedenum}
\makeatother

\begin{document}
\author{K\k{e}stutis \v{C}esnavi\v{c}ius}
\title{Selmer groups as flat cohomology groups}
\date{\today}
\subjclass[2010]{Primary 11G10; Secondary 14F20, 14K02, 14L15}
\keywords{Abelian variety, Selmer group, fppf cohomology, torsor, N\'{e}ron model}
\address{Department of Mathematics, University of California, Berkeley, CA 94720-3840, USA}
\email{kestutis@berkeley.edu}
\urladdr{http://math.berkeley.edu/~kestutis/}

\begin{abstract} Given a prime number $p$, Bloch and Kato showed how the $p^\infty$-Selmer group of an abelian variety $A$ over a number field $K$ is determined by the $p$-adic Tate module. In general, the $p^m$-Selmer group $\Sel_{p^m} A$ need not be determined by the mod $p^m$ Galois representation $A[p^m]$; we show, however, that this is the case if $p$ is large enough. More precisely, we exhibit a finite explicit set of rational primes $\Sigma$ depending on $K$ and $A$, such that $\Sel_{p^m} A$ is determined by $A[p^m]$ for all $p \not \in \Sigma$. In the course of the argument we describe the flat cohomology group $H^1_\fppf(\cO_K, \cA[p^m])$ of the ring of integers of $K$ with coefficients in the $p^m$-torsion $\cA[p^m]$ of the N\'{e}ron model of $A$ by local conditions for $p\not\in \Sigma$, compare them with the local conditions defining $\Sel_{p^m} A$, and prove that $\cA[p^m]$ itself is determined by $A[p^m]$ for such $p$. Our method sharpens the known relationship between $\Sel_{p^m} A$ and  $H^1_\fppf(\cO_K, \cA[p^m])$ and continues to work for other isogenies $\phi$ between abelian varieties over global fields provided that $\deg \phi$ is constrained appropriately. To illustrate it, we exhibit resulting explicit rank predictions for the elliptic curve $11A1$ over certain families of number fields.
\end{abstract} 

\maketitle

\section{Introduction}

Let $K$ be a number field, let $A$ be a $g$-dimensional abelian variety over $K$, and let $p$ be a prime number. Fix a separable closure $K^s$ of $K$. Tate conjectured \cite{Tat66}*{p.~134} that the \emph{$p$-adic Tate module} $T_pA \ce \varprojlim A[p^m](K^s)$ determines $A$ up to an isogeny of degree prime to $p$, and Faltings proved this in \cite{Fal83}*{\S1 b)}%\footnote{By \cite{Tat66}*{Lemmas 1 and 3}, the quoted result by Faltings implies the bijectivity of $\bZ_p \tensor \Hom(A, B) \ra \Hom_{\Gal(\ov{K}/K)}(T_pA, T_pB)$ for all abelian varieties $A, B$ over $K$. In particular, if $\iota\colon T_pA \isomto T_pB$, one can find an isogeny $\phi\colon A\ra B$ whose reduction $\bmod\, p$ agrees with $\iota \bmod p$, hence $p \nmid \deg \phi$.}
. One can ask whether $A[p]$ alone determines $A$ to some extent. Consideration of the case $g = 1$, $p = 2$ shows that for small $p$ one cannot expect much in this direction. However, at least if $g = 1$ and $K = \bQ$, for $p$ large enough (depending on $A$) the Frey--Mazur conjecture \cite{Kra99}*{Conj.~3} predicts that $A[p]$ should determine $A$ up to an isogeny of degree %\footnote{The condition on the degree can be added, because by \cite{Zar85}*{Thm.~1}, up to isomorphism there are only finitely many abelian varieties $B$ isogenous to $A$ over $K$.} 
prime to $p$.

Consider now the $p^\infty$-Selmer group 
\[
\Sel_{p\I} A \subset H^1(K, A[p\I]),
\]
which consists of the classes of cocycles whose restrictions lie in 
\[
A(K_v) \tensor \bQ_p/\bZ_p \subset H^1(K_v, A[p\I])
\]
for every place $v$ of $K$. Note that $A[p^\infty](K^s) = V_pA/T_pA$ with $V_pA \ce T_pA \tensor_{\bZ_p} \bQ_p$, so $T_pA$ determines the Galois cohomology groups appearing in the definition of $\Sel_{p\I} A$. Since an isogeny of degree prime to $p$ induces an isomorphism on $p\I$-Selmer groups, the theorem of Faltings implies that $T_pA$ determines $\Sel_{p\I}A$ up to isomorphism. One may expect, however, a more direct and more explicit description of $\Sel_{p\I}A$ in terms of $T_pA$. For this, it suffices to give definitions of the subgroups $A(K_v) \tensor \bQ_p/\bZ_p \subset H^1(K_v, A[p\I])$ in terms of $T_pA$.

Bloch and Kato found the desired definitions in \cite{BK90}*{\S3}: if $v \nmid p$, then $A(K_v) \tensor \bQ_p/\bZ_p = 0$; if $v \mid p$, then, letting $B_\cris$ be the crystalline period ring of Fontaine and working with Galois cohomology groups formed using continuous cochains in the sense of \cite{Tat76}*{\S2}, they define 
\[ H^1_f(K_v, V_pA) \ce \Ker\p{H^1(K_v, V_pA) \ra H^1(K_v, V_pA \tensor_{\bQ_p} B_\cris)},\]
and prove that 
\[ A(K_v) \tensor \bQ_p/\bZ_p = \im\p{H^1_f(K_v, V_pA) \ra H^1(K_v, V_pA/T_pA) = H^1(K_v, A[p\I])}.\] Considering the $p$-Selmer group $\Sel_pA$ and $A[p]$ instead of $\Sel_{p\I}A$ and $A[p\I]$ (equivalently, $\Sel_{p\I}A$ and $T_pA$), in the light of the Frey--Mazur conjecture, one may expect a direct description of $\Sel_pA$ in terms of $A[p]$ for large $p$. We give such a description as a special case of the following theorem.

%Let $\Sel_{p\I} A \subset H^1(K, A[p\I])$ be the $p^\infty$-Selmer group. The theorem of Faltings implies that $T_pA$ determines $\Sel_{p\I}A$ up to isomorphism, and a direct description of $\Sel_{p\I}A$ in terms of $T_pA$ was given by Bloch and Kato \cite{BK90} using ideas from $p$-adic Hodge theory. Considering the $p$-Selmer group $\Sel_pA$ and $A[p]$ instead of $\Sel_{p\I}A$ and $A[p\I]$ (equivalently, $\Sel_{p\I}A$ and $T_pA$), in light of the Frey--Mazur conjecture, one may expect a direct description of $\Sel_pA$ in terms of $A[p]$ for large $p$. We give such a description as a special case of

\bthm \lab{sel-ind-nf} 
Fix an extension of number fields $L/K$, fix a $K$-isogeny $\phi\colon A \ra B$ between abelian varieties, and let $\cA[\phi]$ and $\cA^L[\phi]$ be the kernels of the induced homomorphisms between the N\'{e}ron models over the rings of integers $\cO_K$ and $\cO_L$. Let $v$ (resp.,~$w$) denote a place of $K$ (resp.,~$L$). For $v, w\nmid \infty$, let $e_v$ and $p_v$ be the absolute ramification index and the residue characteristic of $v$, and let $c_{A, v}$ and $c_{B, v}$ (resp.,~$c_{A, w}$ and $c_{B, w}$) be the local Tamagawa factors of $A$ and $B$.
\benum
\item \lab{sel-ind-nf-aa}
\begin{enumerate}[label={(\roman*)}] 
\item \lab{sel-ind-nf-a} \upshape{(\Cref{sel-cart}, Remark \ref{no-flat}, and \Cref{Aphi-aff}.)}~\it
The pullback map 
\[
H^1_\fppf(\cO_K, \cA[\phi]) \ra H^1(K, A[\phi])
\]
is an isomorphism onto the preimage of $\prod_{v\nmid \infty} H^1_\fppf(\cO_v, \cA[\phi]) \subset \prod_{v\nmid \infty} H^1(K_v, A[\phi])$.

\item \lab{sel-ind-nf-b} \upshape{(\Cref{sel-fppf-gl}~\ref{sel-fppf-gl-e}.)}~\it 
Assume that $A$ has semiabelian reduction at all $v \mid \deg \phi$.
If $\deg \phi$ is prime to $\prod_{v\nmid \infty} c_{A, v}c_{B, v}$ and either $2\nmid \deg \phi$ or $A(K_v)$ is connected for all real $v$, then 
\[
H^1_\fppf(\cO_K, \cA[\phi]) = \Sel_\phi A
\]
inside $H^1(K, A[\phi])$.
\eenum

\item \lab{sel-ind-nf-c} \upshape{(\Cref{Aphi-detd}.)}~\it 
If $A$ has good reduction at all $v \mid \deg \phi$ and if $e_v < p_v - 1$ for every such $v$, then the $\cO_L$-group scheme $\cA^L[\phi]$ is determined up to isomorphism by the $\Gal(L^s/K)$-module $A[\phi](L^s)$.

\eenum
Thus, if $(\deg \phi, \prod_{w\nmid \infty} c_{A, w}c_{B, w}) = 1$, the reduction of $A$ is good at all $v \mid \deg \phi$, and $e_v < p_v - 1$ for every such $v$ (in particular, $2\nmid \deg \phi$), then the $\phi$-Selmer group 
\[
\Sel_\phi A_L \subset H^1(L, A[\phi])
\] 
is determined by the $\Gal(L^s/K)$-module $A[\phi](L^s)$.  \ethm

\bcor
If $A$ has potential good reduction at every finite place of $K$ and $p$ is large enough (depending on $A$), then $A[p^m]$ determines $\Sel_{p^m} A_L$ for every finite extension $L/K$.
\ecor

\bpf
By a theorem of McCallum \cite{ELL96}*{pp.~801--802}, every prime $q$ dividing some $c_{A, w}$ satisfies $q\le 2g +1$. Therefore, it suffices to consider those $p$ with $p > \max(2g + 1, [K : \bQ] + 1)$ for which $A$ has good reduction at every place of $K$ above $p$ and to apply \Cref{sel-ind-nf} to the multiplication by $p^m$ isogeny.
\epf

\brems
\remi
Relationships similar to \ref{sel-ind-nf-b} between Selmer groups and flat cohomology groups are not new and have been implicitly observed already in \cite{Maz72} and subsequently used by Mazur, Schneider, Kato, and others (often after passing to $p\I$-Selmer groups as is customary in Iwasawa theory). However, the description of $H^1_\fppf(\cO_K, \cA[\phi])$ by local conditions in \ref{sel-ind-nf-a} seems not to have appeared in the literature before, and consequently \ref{sel-ind-nf-b} is more precise than what seems to be available.  

In a more restrictive setup, the question of the extent to which $A[\phi]$ determines $\Sel_\phi A$ has also been discussed in \cite{Gre10}.

\remi 
In the case of elliptic curves, Mazur and Rubin find in \cite{MR15}*{Thm.~3.1 and 6.1} (see also \cite{AS05}*{6.6} for a similar result of Cremona and Mazur) that under assumptions different from those of \Cref{sel-ind-nf}, $p^m$-Selmer groups are determined by mod $p^m$ Galois representations together with additional data including the set of places of potential multiplicative reduction. It is unclear to us whether their results can be recovered from the ones presented in this paper.

\remi \lab{sel-coeffs}
The Selmer type description as in \ref{sel-ind-nf-a} continues to hold for $H^1_\et(\cO_K, \cA)$, where $\cA \ra \Spec \cO_K$ is the N\'{e}ron model of $A$. This leads to a reproof of the \'{e}tale cohomological interpretation of the Shafarevich--Tate group $\Sha(A)$ in \Cref{mazur}; such an interpretation is implicit already in the arguments of \cite{Ray65}*{II.\S3} and is proved in \cite{Maz72}*{Appendix}. Our argument seems more direct: in the proof of loc.~cit.~the absence of \Cref{sel-cart} is circumvented with a diagram chase that uses cohomology with supports exact sequences.

\remi
In \Cref{sel-ind-nf}~\ref{sel-ind-nf-aa}, it is possible to relate $\Sel_\phi A$ and $H^1_\fppf(\cO_K, \cA[\phi])$ under weaker hypotheses than those of \ref{sel-ind-nf-b} by combining \Cref{3-sub-comp} with \Cref{sel-cart} as in the proof of \Cref{sel-fppf-gl} (see also Remark \ref{quant}).

\remi
The interpretation of Selmer groups as flat cohomology groups is useful beyond the case when $\phi$ is multiplication by an integer. For an example, see the last sentence of Remark \ref{rem-last}.

\remi 
\Cref{sel-ind-nf} is stronger than its restriction to the case $L = K$. Indeed, the analogue of $e_v < p_v - 1$ may fail for $L$ but hold for $K$. This comes at the expense of $\cA^L[\phi]$ and $\Sel_\phi A_L$ being determined by $A[\phi](L^s)$ as a $\Gal(L^s/K)$-module, rather than as a $\Gal(L^s/L)$-module.

\remi 
Taking $L = K$ and $A = B$ in \Cref{sel-ind-nf}, we get the set $\Sigma$ promised in the abstract by letting it consist of all primes below a place of bad reduction for $A$, all primes dividing a local Tamagawa factor of $A$, the prime $2$, and all odd primes $p$ ramified in $K$ for which $e_v \ge p  - 1$ for some place $v$ of $K$ above $p$.  

\remi In \Cref{sel-ind-nf}, is the subgroup $B(L)/\phi A(L)$ (equivalently, the quotient $\Sha(A_L)[\phi]$) also determined by $A[\phi](L^s)$? The answer is `no'. Indeed, in \cite{CM00}*{p.~24} Cremona and Mazur report\footnote{Cremona and Mazur assume the Birch and Swinnerton-Dyer conjecture to compute Shafarevich--Tate groups analytically. This is unnecessary for us, since full $2$-descent finds provably correct ranks of $2534 E 1$, $2534 G 1$, $4592D1$, and $4592G1$.} that the elliptic curves $2534 E 1$ and $2534 G 1$ over $\bQ$ have isomorphic mod $3$ representations, but $2534 E 1$ has rank $0$, whereas $2534 G 1$ has rank $2$. Since $3$ is prime to the conductor $2534$ and the local Tamagawa factors $c_2 = 44$, $c_7 = 1$, $c_{181} = 2$ (resp., $c_2 = 13$, $c_7 = 2$, $c_{181} = 1$) of $2534 E 1$ (resp.,~$2534 G 1$), \Cref{sel-ind-nf} indeed applies to these curves. Another example (loc.~cit.) is the pair $4592D1$ and $4592G1$ with $\phi = 5$ and ranks $0$ and $2$. %over $\bQ$ that is reported (loc.~cit.)~to have isomorphic mod $5$ representations, ranks $0$ and $2$, and Shafarevich--Tate groups isomorphic to $\bZ/5\bZ \oplus \bZ/5\bZ$ and $0$, respectively. \Cref{sel-ind-nf} applies, because $5$ is prime to the conductor $4592$ and the local Tamagawa factors $c_2 = 4$, $c_7 = 1$, $c_{41} = 2$ (resp., $c_2 = 4$, $c_7 = 1$, $c_{181} = 1$) of $4592 D 1$ (resp.,~$4592 G 1$). 

For an odd prime $p$ and elliptic curves $E$ and $E\pr$ over $\bQ$ with $E[p] \cong E\pr[p]$ and prime to $p$ conductors and local Tamagawa factors, \Cref{sel-ind-nf}, expected finiteness of $\Sha$, and Cassels--Tate pairing predict that $\rk E(\bQ) \equiv \rk E\pr(\bQ) \bmod 2$. Can one prove this directly?

\remi \lab{sel-ind-ff}
For the analogue of \Cref{sel-ind-nf}~\ref{sel-ind-nf-aa} in the case when the base is a global function field, one takes a (connected) proper smooth curve $S$ over a finite field in the references indicated in the statement of \Cref{sel-ind-nf}~\ref{sel-ind-nf-aa}. Letting $K$ be the function field of $S$, the analogue of \Cref{sel-ind-nf}~\ref{sel-ind-nf-c} is \Cref{Aphi-ner}: if $\Char K \nmid \deg \phi$, then $\cA[\phi] \ra S$ is the N\'{e}ron model of $A[\phi] \ra \Spec K$ ($L$ plays no role); in this case, due to \Cref{unr-comp}~\ref{unr-comp-d}, 
\[
H^1_\fppf(S, \cA[\phi]) \subset H^1(K, A[\phi])
\] 
is the subset of the everywhere unramified cohomology classes. The final conclusion becomes: if $(\deg \phi, \Char K \prod_s c_{A, s}c_{B, s}) = 1$ (the product of the local Tamagawa factors is indexed by the closed $s \in S$), then $A[\phi]$ determines the $\phi$-Selmer subgroup
\[
\Sel_\phi A \subset H^1(K, A[\phi]),
\]
which, in fact, consists of the everywhere unramified cohomology classes of~$H^1(K, A[\phi])$.
\erems

\beg
We illustrate our methods and results by estimating the $5$-Selmer group of the base change $E_K$ of the elliptic curve $E = 11 A 1$ to any number field $K$. This curve has also been considered by Tom Fisher, who described in \cite{Fis03a}*{2.1} the $\phi$-Selmer groups of $E_K$ for the two degree $5$ isogenies $\phi$ of $E_K$ defined over $\bQ$. We restrict to $11A1$ for the sake of concreteness (and to get precise conclusions \ref{calc-0}--\ref{calc-e}); our argument leads to estimates analogous to \eqref{11A1} for every elliptic curve $A$ over $\bQ$ and an odd prime $p$ of good reduction for $A$ such that $A[p] \cong \bZ/p\bZ \oplus \mu_p$.

Let $\cE^{K} \ra \Spec \cO_K$ be the N\'{e}ron model of $E_K$. Since $E[5] \cong \bZ/5\bZ \oplus \mu_5$, %Since the local Tamagawa factor $c_{11} = 5$, the separated quasi-finite flat (cf.~\Cref{qf-flat}) $\bZ_{11}$-group scheme $\cE^{\bQ}[5]_{\bZ_{11}}$ has constant fiber degree, hence is finite \cite{DR73}*{II.1.19}. 
 the proof of \Cref{Aphi-detd} supplies an isomorphism
\[
\cE^K[5] \simeq \underline{\bZ/5\bZ}_{\cO_K} \oplus \mu_5.
\] 
Therefore, the cohomology sequence of $0 \ra \mu_5 \ra \bG_m \xra{5} \bG_m \ra 0$ together with the isomorphism $H^1_\fppf(\cO_K, \bZ/5\bZ) \simeq \Cl_K[5]$ gives
\be \lab{11A1-p}
\dim_{\bF_5} H^1_\fppf(\cO_K, \cE^K[5]) = 2\dim_{\bF_5} \Cl_K[5] + \dim_{\bF_5} \cO_K^\times/\cO_K^{\times 5} = 2h_5^K + r_1^K + r_2^K - 1 + u_5^K,
\ee
where $\Cl_K$ is the ideal class group, $r_1^K$ and $r_2^K$ are the numbers of real and complex places, and 
\[
h_5^K \ce \dim_{\bF_5} \Cl_K[5], \quad\quad\quad u_5^K \ce \dim_{\bF_5} \mu_5(\cO_K).
\] 
The component groups of N\'{e}ron models of elliptic curves with split multiplicative reduction are cyclic, so \eqref{11A1-p} and Remark \ref{quant} give the bounds
\be \lab{11A1}
2h_5^K + r_1^K + r_2^K - 1 + u_5^K -
\#\{ v\mid 11\} \le \dim_{\bF_5} \Sel_5 E_K \le 2h_5^K + r_1^K + r_2^K - 1 + u_5^K + \#\{ v\mid 11\}.
\ee
Thus, the obtained estimate is most precise when $K$ has a single place above $11$. Also, 
\be \lab{par}
\dim_{\bF_5} \Sel_5 E_K \equiv r_1^K + r_2^K - 1 + u_5^K + \#\{ v\mid 11\} \bmod 2,
\ee
because the $5$-parity conjecture is known for $E_K$ by the results of \cite{DD08}. When $K$ ranges over the quadratic extensions of $\bQ$, due to \eqref{11A1}, the conjectured unboundedness of the $5$-ranks $h_5^K$ of the ideal class groups is equivalent to the unboundedness of $\dim_{\bF_5} \Sel_5 E_K$. This equivalence is an instance of a general result \cite{Ces15}*{1.5} that gives a precise relation between unboundedness questions for Selmer groups and class groups. That a relation of this sort may be feasible has also been (at least implicitly) observed by other authors, see, for instance, \cite{Sch96}.

It is curious to draw some concrete conclusions from \eqref{11A1} and \eqref{par}.
\benum
\item \lab{calc-0}
As is also well known, $\rk E(\bQ) = 0$.

\item\lab{calc-b} If $K$ is imaginary quadratic with $h_5^K = 0$ and $11$ is inert or ramified in $K$, then $\rk E(K) = 0$.

\item\lab{calc-c} If $K$ is imaginary quadratic with $h_5^K = 0$ and $11$ splits in $K$, then either $\rk E(K) = 1$, or $\rk E(K) = 0$ and $\cork_{\bZ_5}\Sha (E_K)[5\I] = 1$, because, due to the Cassels--Tate pairing, 
\[
\cork_{\bZ_5}\Sha (E_K)[5\I]  \equiv \dim_{\bF_5} \Sha(E_K)[5] \bmod 2.
\]
Mazur in \cite{Maz79}*{Thm.~on p.~237} and Gross in \cite{Gro82}*{Prop.~3} proved that $\rk E(K) = 1$.

\item \lab{calc-d} If $F$ is a quadratic extension of a $K$ as in \ref{calc-c} in which none of the places of $K$ above $11$ split and $h_5^F = 0$, then either $\rk E(F) = 2$, or $\Sha (E_F)[5\I]$ is infinite (one again uses the Cassels--Tate pairing).

\item \lab{calc-real} If $K$ is real quadratic with $h_5^K = 0$ and $11$ is inert or ramified in $K$, then either $\rk E(K) = 1$, or $\rk E(K) = 0$ and $\cork_{\bZ_5} \Sha (E_K)[5\I] = 1$. In the latter case $\Sha(E_K)[p\I]$ is infinite for every prime $p$, because the $p$-parity conjecture is known for $E_K$ for every $p$ by \cite{DD10}*{1.4} (applied to $E$ and its quadratic twist by $K$). Gross proved in \cite{Gro82}*{Prop.~2} that if $11$ is inert, then $\rk E(K) = 1$.

\item\lab{calc-e} If $K$ is cubic with a complex place (or quartic totally imaginary), a single place above $11$, and $h_5^K = 0$, then either $\rk E(K) = 1$, or $\rk E(K) = 0$ and $\cork_{\bZ_5} \Sha (E_K)[5\I] = 1$.
\eenum

How can one construct the predicted rational points? In \ref{calc-c} and the inert case of \ref{calc-real}, \cite{Gro82} explains that Heegner point constructions account for the predicted rank growth. %However, \ref{calc-d} and \ref{calc-e} concern situations that seem to be beyond the scope of applicability of the existing methods for systematic construction of rational points of infinite order.
\eeg

\bpp[The contents of the paper] 
We begin by restricting to local bases in \S\ref{local} and comparing the subgroups $B(K_v)/\phi(A(K_v))$, $H^1_\fppf(\cO_v, \cA[\phi])$, and $H^1_\nr(K_v, A[\phi])$ of $H^1(K_v, A[\phi])$ under appropriate hypotheses. In \S\ref{glue}, after recording some standard results on fpqc descent, we apply them to prove \Cref{sel-ind-nf}~\ref{sel-ind-nf-c} and to reprove the \'{e}tale cohomological interpretation of Shafarevich--Tate groups. In \S\ref{sel-type}, exploiting the descent results of \S\ref{glue}, we take up the question of $H^1_\fppf$ with appropriate coefficients over Dedekind bases being described by local conditions and prove \Cref{sel-ind-nf}~\ref{sel-ind-nf-a}. The final \S\ref{selmer-flat} uses the local analysis of \S\ref{local} to compare $\Sel_\phi A$ and $H^1_\fppf(\cO_K, \cA[\phi])$ and to complete the proof of \Cref{sel-ind-nf}. The two appendices collect various results concerning torsors and exact sequences of N\'{e}ron models used in the main body of the text. 

Some of the results presented in this paper are worked out in somewhat more general settings in the PhD thesis of the author; we invite a reader interested in this to consult \cite{Ces14}, which also discusses several tangentially related questions.
\epp

\bpp[Conventions]\lab{conv} When needed, a choice of a separable closure $K^s$ of a field $K$ will be made implicitly, as will be a choice of an embedding $K^s \hra L^s$ for an overfield $L/K$. If $v$ is a place of a global field $K$, then $K_v$ is the corresponding completion; for $v\nmid \infty$, the ring of integers and the residue field of $K_v$ are denoted by $\cO_v$ and $\bF_v$. If $K$ is a number field, $\cO_K$ is its ring of integers. For $s\in S$ with $S$ a scheme, $\cO_{S, s}$, $\fm_{S, s}$, and $k(s)$ are the local ring at $s$, its maximal ideal, and its residue field. For a local ring $R$, its henselization, strict henselization, and completion are $R^h$, $R^{sh}$, and $\wh{R}$. The fppf, big \'{e}tale, and \'{e}tale sites of $S$ are $S_\fppf$, $S_\Et$, and $S_\et$; the objects of $S_\fppf$ and $S_\Et$ are all $S$-schemes, while those of $S_\et$ are all schemes \'{e}tale over $S$. The cohomology groups computed in $S_\et$ and $S_\fppf$ are denoted by $H^i_\et(S, -)$ and $H^i_\fppf(S, -)$; Galois cohomology merits no subscript: $H^i(K, -)$. An fppf torsor is a torsor under the group in question for the fppf topology. An algebraic group over a field $K$ is a smooth $K$-group scheme of finite type. 
\epp

\subsection*{Acknowledgements} I thank Bjorn Poonen for many helpful discussions, suggestions, and for reading various drafts. I thank Brian Conrad for reading the manuscript and suggesting numerous improvements. I thank the referee for helpful suggestions that improved the manuscript. I thank Rebecca Bellovin, Henri Darmon, Tim Dokchitser, Jessica Fintzen, Jean Gillibert, Benedict Gross, Mark Kisin, Chao Li, Dino Lorenzini, Barry Mazur, Martin Olsson, Michael Stoll, and David Zureick-Brown for helpful conversations or correspondence regarding the material of this paper. Part of the research presented here was carried out during the author's stay at the Centre Interfacultaire Bernoulli (CIB) in Lausanne during the course of the program ``Rational points and algebraic cycles''. I thank CIB, NSF, and the organizers of the program for a lively semester and the opportunity to take part.

\section{Images of local Kummer homomorphisms as flat cohomology groups}\lab{local}

Let $S = \Spec \fo$ for a Henselian discrete valuation ring $\fo$ with a finite residue field $\bF$, let $k = \Frac \fo$, let $i\colon \Spec \bF \ra S$ be the closed point, let $\phi\colon A \ra B$ be a $k$-isogeny of abelian varieties, let $\phi\colon \cA \ra \cB$ be the induced $S$-homomorphism between the N\'{e}ron models, which gives rise to the homomorphism $\phi\colon \Phi_A \ra \Phi_B$ between the \'{e}tale $\bF$-group schemes of connected components of $\cA_\bF$ and $\cB_\bF$. We use various open subgroup schemes of $\cA$ and $\cB$ discussed in \S\ref{open}.

\bpp[The three subgroups]\lab{3-subgp} The first subgroup of $H^1_\fppf(k, A[\phi])$ is 
\[
B(k)/\phi A(k) \cong \im\p{B(k) \xra{\kappa_\phi} H^1_\fppf(k, A[\phi])} \subset H^1_\fppf(k, A[\phi]).
\]

The second subgroup is 
\[
H^1_\fppf(\fo, \cA[\phi])\cong \im\p{H^1_\fppf(\fo, \cA[\phi]) \xra{a} H^1_\fppf(k, A[\phi])} \subset H^1_\fppf(k, A[\phi]),
\]
where the isomorphism results from the injectivity of $a$ supplied by \Cref{Aphi-aff}, \cite{GMB13}*{Prop.~3.1}, and \Cref{sub-ner-pf} (even though $\cA[\phi]$ may fail to be flat, loc.~cit.~proves that its category of fppf torsors is equivalent to the category of fppf torsors of the $\fo$-flat schematic image of $A[\phi]$ in $\cA$, so \Cref{sub-ner-pf} nevertheless applies). 

The third is the unramified subgroup 
\[
H^1_\nr(k, A[\phi]) \ce \Ker\p{H^1(k, A[\phi]) \ra H^1(k^{sh}, A[\phi])} \subset H^1(k, A[\phi]),
\]
where $k^{sh} \ce \Frac \fo^{sh}$. The unramified subgroup is of most interest in the case when $A[\phi]$ is \'{e}tale (for instance, when $\Char k \nmid \deg \phi$); beyond this \'{e}tale case, the unramified subgroup is often too small in comparison to the first two subgroups. 

While $\im \kappa_\phi$ is used to define the $\phi$-Selmer group, $H^1_\fppf(\fo, \cA[\phi])$ and $H^1_\nr(k, A[\phi])$ are easier to study because they depend only on $\cA[\phi]$. We investigate $\im \kappa_\phi$ by detailing its relations with $H^1_\fppf(\fo, \cA[\phi])$ and $H^1_\nr(k, A[\phi])$ in \Cref{3-sub-comp,unr-comp}.  \epp

\blem \lab{Lang-thm} For a commutative connected algebraic group $G \ra \Spec \bF$, one has 
\[
H^j(\bF, G) = 0 \quad\quad \text{for $j \ge 1$.}
\] 
\elem

\bpf 
In the case $j = 1$, the claimed vanishing is a well-known result of Lang \cite{Lan56}*{Thm.~2}. In the case $j > 1$, the vanishing follows from the facts that $\bF$ has cohomological dimension $1$ and that $G(\bF^s)$ is a torsion group (the latter results from the finiteness of $\bF$).
\epf

\blem\lab{nat-iso} For an $\bF$-subgroup $\Gamma \subset \Phi_A$, pullback induces isomorphisms 
\[
H^j_\fppf(\fo, \cA^\Gamma) \cong H^j(\bF, \Gamma) \quad \quad \text{for $j \ge 1$.}
\] 
In particular, $\#H^1_\fppf(\fo, \cA^\Gamma) = \#\Gamma(\bF)$ and $H^j_\fppf(\fo, \cA^{\Gamma}) = 0$ for $j \ge 2$.  \elem

\bpf 
By \cite{Gro68}*{11.7 2$^\circ$)}, pullback induces isomorphisms 
\[
H^j_\fppf(\fo, \cA^\Gamma) \cong H^j(\bF, \cA^{\Gamma}_\bF) \quad \quad \text{ for $j \ge 1$,} 
\]
so it remains to apply \Cref{Lang-thm} to the terms $H^j(\bF, \cA^0_\bF)$ in the long exact cohomology sequence of 
\[
0 \ra \cA^0_\bF \ra \cA^\Gamma_\bF \ra \Gamma \ra 0. \qedhere
\]
\epf

\bpp[The local Tamagawa factors]\lab{loc-tam} These are 
\[
c_A \ce \#\Phi_A(\bF) \quad\quad \text{and} \quad\quad c_B \ce \#\Phi_B(\bF).
\]
The sequences
\[\ba 
0 \ra \Phi_A[\phi](\bF^s) \ra &\Phi_A(\bF^s) \ra (\phi(\Phi_A))(\bF^s) \ra 0, \\
 0 \ra (\phi(\Phi_A))(\bF^s) \ra &\Phi_B(\bF^s) \ra (\Phi_B/\phi(\Phi_A))(\bF^s) \ra 0
 \ea\]
are exact, so
\be\lab{tam-bds} \f{\#\Phi_A(\bF)}{\#(\phi(\Phi_A))(\bF)} \le \# \Phi_A[\phi](\bF) \quad\quad \text{and} \quad \quad \f{\#\Phi_B(\bF)}{\#(\phi(\Phi_A))(\bF)} \le \#\p{\f{\Phi_B}{\phi(\Phi_A)}}(\bF).\ee
\epp

We now compare the subgroups $\im \kappa_\phi$ and $H^1_\fppf(\fo, \cA[\phi])$ of $H^1_\fppf(k, A[\phi])$ discussed in \S\ref{3-subgp}.

\bprop\lab{3-sub-comp}
Suppose that $\cA \xra{\phi} \cB$ is flat (e.g., that $\Char \bF \nmid \deg \phi$ or that $A$ has semiabelian reduction, see \Cref{surj-id}).
\benum 
\item \lab{3-sub-comp-a}
Then
\[\ba 
\#\p{\f{H^1_\fppf(\fo, \cA[\phi])}{H^1_\fppf(\fo, \cA[\phi])\cap \im \kappa_\phi} } &=  \f{\#\Phi_A(\bF)}{\#(\phi(\Phi_A))(\bF)} \overset{\eqref{tam-bds}}{\le} \# \Phi_A[\phi](\bF),\\
\#\p{\f{\im \kappa_\phi}{H^1_\fppf(\fo, \cA[\phi])\cap \im \kappa_\phi}  } &= \f{\#\Phi_B(\bF)}{\#(\phi(\Phi_A))(\bF)} \overset{\eqref{tam-bds}}{\le} \#\p{\f{\Phi_B}{\phi(\Phi_A)}}(\bF). 
\ea\]

\item \lab{3-sub-comp-d}  
If $\deg \phi$ is prime to $c_B$, then $\Phi_B(\bF) = (\phi(\Phi_A))(\bF)$, and hence, by \ref{3-sub-comp-a}, 
\[
\im \kappa_\phi \subset H^1_\fppf(\fo, \cA[\phi]).
\]

\item \lab{3-sub-comp-c} 
If $\deg \phi$ is prime to $c_A$, then $\Phi_A(\bF) = (\phi(\Phi_A))(\bF)$, and hence, by \ref{3-sub-comp-a}, 
\[
H^1_\fppf(\fo, \cA[\phi]) \subset \im \kappa_\phi.
\]

\item \lab{3-sub-comp-e}  If  $\deg \phi$ is prime to $c_Ac_B$, then 
\[
\im \kappa_\phi = H^1_\fppf(\fo, \cA[\phi]).
\]

\eenum
\eprop

\bpf \hfill \benum

\item 
The short exact sequence 
\[
0 \ra \cA[\phi] \ra \cA \xra{\phi} \cB^{\phi(\Phi_A)} \ra 0 
\]
of \Cref{ses-mult} gives 
\[
\xymatrix@C=12pt{ 
\quad 0 \ar[r] & \cB^{\phi(\Phi_A)}(\fo)/\phi\cA(\fo) \ar[r]\ar@{^(->}[d]\ar[r] & H^1_\fppf(\fo, \cA[\phi]) \ar[r]\ar[r]\ar@{^(->}[d]^{a} & \Ker\p{H^1_\fppf(\fo, \cA) \xra{H^1_\fppf(\phi)} H^1_\fppf(\fo, \cB^{\phi(\Phi_A)})} \ar[r]\ar@{^(->}[d] & 0 \\
\quad  0 \ar[r] & B(k)/\phi A(k) \ar[r]^-{\kappa_\phi} & H^1_\fppf(k, A[\phi]) \ar[r] & H^1_\fppf(k, A)[\phi] \ar[r] & 0, } 
 \]
where the injectivity of the vertical arrows follows from the N\'{e}ron property, the snake lemma, and \Cref{inj-ner}. By \Cref{nat-iso}, $H^1_\fppf(\phi)$ identifies with 
\[
H^1(\bF, \Phi_A) \xra{h} H^1(\bF, \phi(\Phi_A))
\]
induced by $\phi$; moreover, $h$ is onto. Since 
\[
\f{H^1_\fppf(\fo, \cA[\phi])}{H^1_\fppf(\fo, \cA[\phi])\cap \im \kappa_\phi} \cong \Ker H^1_\fppf(\phi) \cong \Ker h
\]
and 
\[
\#\Ker h =  \f{\#H^1(\bF, \Phi_A)}{\#H^1(\bF, \phi(\Phi_A))} = \f{\#\Phi_A(\bF)}{\#(\phi(\Phi_A))(\bF)},
\] 
the first claimed equality follows. 

On the other hand,
\be\lab{seq1} \f{\im \kappa_\phi}{H^1_\fppf(\fo, \cA[\phi])\cap \im \kappa_\phi} \cong \f{B(k)/\phi A(k)}{\cB^{\phi(\Phi_A)}(\fo)/\phi\cA(\fo)} \cong  \f{\cB(\fo)}{\cB^{\phi(\Phi_A)}(\fo)}. \ee
Moreover, \Cref{nat-iso} and the \'{e}tale cohomology sequence of the short exact sequence
\[
0 \ra \cB^{\phi(\Phi_A)} \ra \cB \ra i_*(\Phi_B/\phi(\Phi_A)) \ra 0
\] 
from \Cref{ses-gam-gam} give the exact sequence (see \cite{Gro68}*{11.7 1$^\circ$)} for the identifications between different cohomology theories)
\be\lab{seq2} 0 \ra \f{\cB(\fo)}{\cB^{\phi(\Phi_A)}(\fo)} \ra \p{\f{\Phi_B}{\phi(\Phi_A)}}(\bF) \ra H^1(\bF, \phi(\Phi_A)) \ra H^1(\bF, \Phi_B) \surjects H^1\p{\bF, \f{\Phi_B}{\phi(\Phi_A)}}, \ee
where we have used the exactness of $i_*$ for the \'{e}tale topology to obtain the last term.
By combining \eqref{seq1} and \eqref{seq2}, we obtain the remaining equality
\[ \quad \quad \quad \#\p{\f{\im \kappa_\phi}{H^1_\fppf(\fo, \cA[\phi])\cap \im \kappa_\phi}  } = \f{\#(\Phi_B/\phi(\Phi_A))(\bF) \cdot \# H^1(\bF, \Phi_B)}{\# H^1(\bF, \phi(\Phi_A)) \cdot \#H^1(\bF, \Phi_B/\phi(\Phi_A))} = \f{\#\Phi_B(\bF)}{\#(\phi(\Phi_A))(\bF)}. \]

\item Let $\psi\colon B \ra A$ be the isogeny with $\ker \psi = \phi(A[\deg \phi])$, so 
\[
\psi \circ \phi = \deg \phi, \quad \quad \text{and thus also} \quad \quad \phi \circ \psi = \deg \phi.
\] 
If $(\deg \phi, \#\Phi_B(\bF)) = 1$, then 
\[
\quad \quad \Phi_B(\bF) = (\deg \phi)(\Phi_B(\bF)) \subset ((\deg \phi)(\Phi_B))(\bF) \subset (\phi(\Phi_A))(\bF) \subset \Phi_B(\bF),
\] 
which gives the desired equality $\Phi_B(\bF)  = (\phi(\Phi_A))(\bF)$.

\item 
We have the inclusion 
\[
\Phi_A[\phi] \subset \Phi_A[\deg \phi],
\] 
so if $(\deg \phi, \#\Phi_A(\bF)) = 1$, then $\Phi_A[\phi](\bF) = 0$. The resulting injection 
\[
\Phi_A(\bF) \hra \phi(\Phi_A)(\bF) 
\] 
is then surjective because $\#H^1(\bF, \Phi_A[\phi]) = \#\Phi_A[\phi](\bF)$ due to the finiteness of $\bF$.

\item The claim follows by combining \ref{3-sub-comp-d} and \ref{3-sub-comp-c}. \qedhere
\eenum \epf

\brem
In the case $\dim A = 1$ and $\phi = p^m$, \Cref{3-sub-comp}~\ref{3-sub-comp-e} has also been observed by Mazur and Rubin in \cite{MR15}*{Prop.~5.8}.
\erem

We now compare the unramified subgroup $H^1_\nr(k, A[\phi]) \subset H^1(k, A[\phi])$ to $\im \kappa_\phi$ and $H^1_\fppf(\fo, \cA[\phi])$:

\bprop\lab{unr-comp} 
Suppose that $A[\phi]$ is \'{e}tale (e.g.,~that $\Char k \nmid \deg \phi$), and let $\cG \ra S$ be the N\'{e}ron model of $A[\phi] \ra \Spec K$ (the N\'{e}ron model exists by, for instance, \cite{BLR90}*{\S7.1, Cor.~6}).
\benum
\item\lab{unr-comp-c} 
There is an inclusion 
\[
\quad H^1_\nr(k, A[\phi]) \subset H^1_\fppf(\fo, \cA[\phi])
\]
inside $H^1(k, A[\phi])$.

\item\lab{unr-comp-d} If $\cA[\phi] \ra S$ is \'{e}tale (e.g.,~if $\Char \bF \nmid \deg \phi$), then 
\[
H^1_\nr(k, A[\phi]) = H^1_\fppf(\fo, \cA[\phi])
\]
inside $H^1(k, A[\phi])$.

\item \lab{unr-comp-f} One has 
\[
H^1_\nr(k, A[\phi]) \subset \im \kappa_\phi
\] 
inside $H^1(k, A[\phi])$ if one assumes in addition that 
\begin{enumerate}[label={(\roman*)}]
\item \lab{ass-i}
$\cA \xra{\phi} \cB$ is flat (which holds if $\Char \bF \nmid \deg \phi$ or if $A$ has semiabelian reduction, see \Cref{surj-id}), and
\item \lab{ass-ii}
% and that %the local Tamagawa factor 
$\#\Phi_A(\bF) = \#(\phi(\Phi_A))(\bF)$ (which holds if $\deg \phi$ is prime to $c_A$, see \Cref{3-sub-comp}~\ref{3-sub-comp-c}). 
\eenum

\item \lab{unr-comp-e} One has 
\[
\quad H^1_\nr(k, A[\phi]) = \im \kappa_\phi = H^1_\fppf(\fo, \cA[\phi])
\]
inside $H^1(k, A[\phi])$ if one assumes in addition that 
\begin{enumerate}[label={(\roman*)}]
\item $\cA[\phi] \ra S$ is \'{e}tale (which holds if $\Char \bF \nmid \deg \phi$), and %the local Tamagawa factors 
\item $\#\Phi_A(\bF) = \#(\phi(\Phi_A))(\bF) = \#\Phi_B(\bF)$ (which holds if $\deg \phi$ is prime to $c_Ac_B$).
\eenum
\eenum
 \eprop

\bpf \hfill
\benum
\item By \Cref{im-ner-unr} (together with \cite{Gro68}*{11.7 1$^\circ$)} for the identification between the \'{e}tale and the fppf cohomology groups), 
\[
\quad H^1_\nr(k, A[\phi]) = H^1_\fppf(\fo, \cG)
\] 
inside $H^1(k, A[\phi])$. It therefore suffices to find an $S$-homomorphism $\cG \ra \cA[\phi]$ that induces an isomorphism on the generic fibers. Such an $S$-homomorphism is provided by \cite{BLR90}*{\S7.1, Cor.~6}, which describes $\cG$ as the group smoothening of the schematic image of $A[\phi]$ in $\cA$.

\item 
If $\cA[\phi] \ra S$ is \'{e}tale, then no smoothening is needed in loc.~cit.,~that is, $\cA[\phi]$ itself is the N\'{e}ron model of $A[\phi]$. The claim therefore results from \Cref{im-ner-unr}.

\item 
Due to the assumptions \ref{ass-i} and \ref{ass-ii}, \Cref{3-sub-comp}~\ref{3-sub-comp-a} applies and gives the inclusion 
\[
H^1_\fppf(\fo, \cA[\phi]) \subset \im \kappa_\phi.
\]
The claim therefore results from \ref{unr-comp-c}.

\item \Cref{3-sub-comp} supplies the equality 
\[
H^1_\fppf(\fo, \cA[\phi]) = \im \kappa_\phi,
\]
so the claim results from \ref{unr-comp-d}. \qedhere
\eenum \epf

\brem \lab{cas-lem}
\Cref{unr-comp}~\ref{unr-comp-e} generalizes a well-known lemma of Cassels \cite{Cas65}*{4.1}, which yields $\im \kappa_\phi = H^1_\nr(k, A[\phi])$ under the assumptions that $\Char \bF \nmid \deg \phi$ and that the reduction is good (so that $c_A = c_B = 1$). In a setting where $\dim A = 1$ and $\Char \bF \nmid \deg \phi$, a special case of this generalization has also been observed by Schaefer and Stoll \cite{SS04}*{proof of Prop.~3.2}. 
\erem

\section{Assembling $\cA[\phi]$ by glueing}\lab{glue}

A descent lemma \ref{eq-cat} formalizes the idea that giving a scheme over a connected Dedekind scheme $S$ amounts to giving a scheme over a nonempty open $V \subset S$ together with a compatible $\wh{O}_{S, s}$-scheme for every $s \in S - V$. \Cref{eq-cat} is crucial for glueing $\cA[\phi]$ together in the proof of \Cref{Aphi-detd}; it will also be key for Selmer type descriptions of sets of torsors in \S\ref{sel-type}. Its more technical part \ref{eq-cat-b} involving algebraic spaces is needed in order to avoid a quasi-affineness hypothesis in \Cref{sel-cart}. This corollary enables us to glue torsors under a N\'{e}ron model in the proof of \Cref{mazur}: even though \emph{a posteriori} such torsors are schemes, we glue them as algebraic spaces because the description of the essential image in \Cref{eq-cat}~\ref{eq-cat-a} is not practical beyond the quasi-affine case. For the proof of \Cref{sel-ind-nf}, however, there is no need to resort to algebraic spaces: \Cref{eq-cat}~\ref{eq-cat-a} is sufficient due to the affineness of $\cA[\phi]$ guaranteed by \Cref{Aphi-aff}.

\blem\lab{eq-cat} Let $R$ be a discrete valuation ring, set $K \ce \Frac R$ and $K^h \ce \Frac R^h$, and consider
\[ F\colon X \mapsto (X_{K}, X_{R^h}, \tau\colon (X_{K})_{K^h} \xra{\sim} (X_{R^h})_{K^h}), \]
a functor from the category of $R$-algebraic spaces to the category of triples consisting of a $K$-algebraic space, an $R^h$-algebraic space, and an isomorphism between their base changes to $K^h$.
\benum
\item\lab{eq-cat-a} When restricted to the full subcategory of $R$-schemes, $F$ is an equivalence onto the full subcategory of triples of schemes that admit a quasi-affine open covering (see the proof for the definition). The same conclusion holds with $R^h$ and $K^h$ replaced by $\wh{R}$ and $\wh{K} \ce \Frac \wh{R}$.

\item\lab{eq-cat-b} When restricted to the full subcategory of $R$-algebraic spaces of finite presentation, $F$ is an equivalence onto the full subcategory of triples involving only algebraic spaces of finite presentation. %If $\wh{K}/K$ is separable\footnote{A field extension $l/k$ is \emph{separable} if for every field extension $L/k$, the ring $L\tensor_k l$ is reduced.} (equivalently, if $R$ is excellent), the same conclusion holds with $R^h$ and $K^h$ replaced by $\wh{R}$ and $\wh{K}$.
\eenum
\elem
   
\bpf \hfill
\benum
\item 
This is proved in \cite{BLR90}*{\S6.2,~Prop.~D.4~(b)}. A triple of schemes admits a quasi-affine open covering if 
\[
\textstyle X_{K} = \bigcup_{i \in I} U_i \quad \quad \text{and}\quad \quad X_{R^h} = \bigcup_{i \in I} V_i
\]
for quasi-affine open subschemes $U_i \subset X_K$ and $V_i \subset X_{R^h}$ for which $\tau$ restricts to isomorphisms $(U_i)_{K^h} \xra{\sim} (V_i)_{K^h}$.

\item The method of proof was suggested to me by Brian Conrad. By construction, $R^h$ is a filtered direct limit of local \'{e}tale $R$-algebras $R\pr$ which are discrete valuation rings sharing the residue field and a uniformizer with $R$. 
%Here's the argument: $\Spec R$ is connected, also regularity ascends and dimension is preserved, so $R$ is a regular local ring of dimension $1$, hence a DVR; sharing a uniformizer comes from \'{e}tale being in particular unramified; sharing residue field comes from the definition of henselization.
Given a 
\[
T = (Y, \cY, \tau\colon Y_{K^h} \isomto \cY_{K^h})
\]
with $Y \ra \Spec K$ and $\cY \ra \Spec R^h$ of finite presentation, to show that it is in the essential image of the restricted $F$, we first use limit considerations (for instance, as in \cite{Ols06}*{proof of Prop.~2.2}) to descend $\cY$ to a $\cY\pr \ra \Spec R\pr$ for some $R\pr$ as above. 

Similarly, $K^h = \varinjlim K\pr$ with $K\pr \ce \Frac R\pr$, so $\tau$ descends to a $\tau\pr\colon Y_{K\pr} \xra{\sim} \cY\pr_{K\pr}$ after possibly enlarging $R\pr$. We transport the $K\pr/K$-descent datum on $Y_{K\pr}$ along $\tau\pr$ to get a descent datum on $\cY\pr_{K\pr}$, which, as explained in \cite{BLR90}*{\S6.2,~proof of Lemma C.2}, extends uniquely to an $R\pr/R$-descent datum on $\cY\pr$. By \cite{LMB00}*{1.6.4}, this descent datum is effective, and we get a quasi-separated $R$-algebraic space $X$; by construction, $F(X) \cong T$, and by \cite{SP}*{\href{http://stacks.math.columbia.edu/tag/041V}{041V}}, $X$ is of finite presentation. 

The full faithfulness of $F$ follows~from a similar limit argument that uses \'{e}tale descent for morphisms of sheaves on $R_\Et$ and \cite{LMB00}*{4.18~(i)}. \qedhere
 \eenum \epf

Let $S$ be a connected Dedekind scheme (see \S\ref{ded-sch} for the definition), let $K$ be its function field. For $s \in S$, set $K_{S, s} \ce \Frac \cO_{S, s}$. The purpose of this convention (note that $K_{S, s} = K$) is to clarify the statement of \Cref{eq-qa-gp} by making $\cO_{S, s}$ and $K_{S, s}$ notationally analogous to $\cO_{S, s}^h$ and $K_{S, s}^h$.

\bcor \lab{eq-qa-gp}
Let $S$ be a Dedekind scheme, let $s_1, \dotsc, s_n \in S$ be distinct nongeneric points, and let $V \ce S - \{ s_1, \dotsc, s_n\}$ be the complementary open subscheme. The functor  
\[
F\colon \cG \mapsto (\cG_V, \cG_{\cO_{S, s_1}}, \dotsc, \cG_{\cO_{S, s_n}}, \gA_i\colon (\cG_V)_{K_{S, s_i}} \isomto (\cG_{\cO_{S, s_i}})_{K_{S, s_i}}\text{ for } 1 \le i \le n)
\]
is an equivalence of categories from the category of quasi-affine $S$-group schemes to the category of tuples consisting of a quasi-affine $V$-group scheme, a quasi-affine $\cO_{S, s_i}$-group scheme for each $i$, and isomorphisms $\gA_1, \dotsc, \gA_n$ of base changed group schemes as indicated. The same conclusion holds with $\cO_{S, s_i}$ and $K_{S, s_i}$ replaced by $\cO_{S, s_i}^h$ and $K_{S, s_i}^h$ or by $\wh{\cO}_{S, s_i}$ and $\wh{K}_{S, s_i}$. 
\ecor

\bpf 
For localizations, the claim is a special case of fpqc descent.
For henselizations and completions, the claim therefore follows from \Cref{eq-cat}.
\epf

\bprop[\Cref{sel-ind-nf}~\ref{sel-ind-nf-c}] \lab{Aphi-detd}
Let $L/K$ be an extension of number fields, and let $\phi\colon A \ra B$ be a $K$-isogeny between abelian varieties. Assume that
\begin{enumerate}[label={(\roman*)}]
\item \lab{Aphi-detd-i} $A$ has good reduction at all the places $v\mid \deg\phi$ of $K$;

\item\lab{Aphi-detd-ii} For every place $v \mid \deg \phi$ of $K$, its absolute ramification index $e_v$ satisfies 
\[
e_v < p_v -1,
\]
where $p_v$ is the residue characteristic of $v$.
\eenum
Then the $\cO_L$-group scheme $\cA^L[\phi]$, defined as the kernel of the homomorphism induced by $\phi_L$ between the N\'{e}ron models over $\cO_L$, is determined up to isomorphism by the $\Gal(L^s/K)$-module $A[\phi](L^s)$.
\eprop

\bpf 
By \Cref{Aphi-ner}, $\cA^L[\phi]_{S[\f{1}{\deg \phi}]}$ is the N\'{e}ron model of the finite \'{e}tale $A[\phi]_L$, and hence is determined by $A[\phi]$. By \Cref{eq-qa-gp}, it therefore suffices to prove that each $\cA^L[\phi]_{\cO_w}$ for a place $w\mid \deg \phi$ of $L$ is also determined by $A[\phi]$. Moreover, if such a $w$ lies above the place $v$ of $K$, then the good reduction assumption implies that 
\[
\cA^L[\phi]_{\cO_w} \cong (\cA^K[\phi]_{\cO_v})_{\cO_w},
\]
so it suffices to prove that already $\cA^K[\phi]_{\cO_v}$ is determined by $A[\phi]$.

Let $p$ be the residue characteristic of $v$. By \Cref{qf-flat}, $\cA^K[\phi]_{\cO_v}$ is finite flat, so it uniquely decomposes as a direct product of commutative finite flat $\cO_v$-group schemes of prime power order. The prime-to-$p$ factor is finite \'{e}tale, so it is the N\'{e}ron model of the prime-to-$p$ factor of $A[\phi]$, and hence is determined by $A[\phi]$. The $p$-primary factor is also determined thanks to Raynaud's result \cite{Ray74}*{Thm.~3.3.3} on uniqueness of finite flat models over Henselian discrete valuation rings of mixed characteristic and low absolute ramification index.
\epf

\brem \lab{wo-ii} Dropping \ref{Aphi-detd-ii} but keeping \ref{Aphi-detd-i}, the proof continues to give the same conclusion as long as one argues that in the situation at hand $\cA^K[\phi]_{\cO_v}$ is determined by $A[\phi]$ for each $v \mid \deg \phi$.
\erem

Although the assumption \ref{Aphi-detd-ii} excludes the cases when $2 \mid \deg \phi$, Remark \ref{wo-ii} can sometimes be used to overcome this, as the following example illustrates.

\beg \lab{2-ss-eg} 
Let $K$ be a number field of odd discriminant, and let $A \ra \Spec K$ be an elliptic curve with good reduction at all $v\mid 2$. Assume that $A[2](K_v) \neq (\bZ/2\bZ)^2$ for every $v \mid 2$, so that $A[2]_{K_v}$ has at most one $K_v$-subgroup of order $2$ for every such $v$.
We show that the conclusion of \Cref{Aphi-detd} holds for $2\colon A \ra A$, so, in particular, if $\prod_{v\nmid \infty} c_{A, v}$ is odd and $K$ is totally imaginary, then $A[2]$ determines the $2$-Selmer group $\Sel_2 A$ by \Cref{sel-ind-nf}. 

\Cref{wo-ii} reduces us to proving that $A[2]_{K_v}$ determines $\cA^K[2]_{\cO_v}$ for each $v\mid 2$. We analyze the ordinary and the supersingular reduction cases separately. This is permissible because these cases are distinguishable: in the former, $A[2]_{K_v}$ is reducible, whereas in the latter it is not.

In the supersingular case, by \cite{Ser72}*{p.~275, Prop.~12}, $A[2]_{K_v^{sh}}$ with $K_v^{sh} \ce \Frac \cO_v^{sh}$ is irreducible and also an $\bF_4$-vector space scheme of dimension $1$. By \cite{Ray74}*{3.3.2~3$^{\text{o}}$}, $\cA^K[2]_{\cO_v^{sh}}$ is its unique finite flat $\cO_v^{sh}$-model. By schematic density considerations, the descent datum on $\cA^K[2]_{\cO_v^{sh}}$ with respect to $\cO_v^{sh}/\cO_v$ is uniquely determined by its restriction to the generic fiber, which in turn is determined by $A[2]_{K_v}$. Fpqc descent along $\cO_v^{sh}/\cO_v$ then implies that $A[2]_{K_v}$ determines $\cA^K[2]_{\cO_v}$.

In the ordinary case, the connected-\'{e}tale decomposition shows that $\cA^K[2]_{\cO_v}$ is an extension of $\underline{\bZ/2\bZ}_{\cO_v}$ by $(\mu_2)_{\cO_v}$. Therefore, since we assumed that $A[2]_{K_v}$ determines its subgroup $(\mu_2)_{K_v}$, it also determines $\cA^K[2]_{\cO_v}$ due to the injectivity of 
\[
\Ext_{\cO_v}^1(\bZ/2\bZ, \mu_2) \cong H^1_\fppf(\cO_v, \mu_2) \ra H^1_\fppf(K_v, \mu_2) \cong \Ext_{K_v}^1(\bZ/2\bZ, \mu_2)
\]
(extensions in the category of fppf sheaves of $\bZ/2\bZ$-modules).
\eeg

\section{Selmer type descriptions of sets of torsors}\lab{sel-type}

The main result of this section is \Cref{sel-cart}, which describes certain sets of torsors by local conditions and proves \Cref{sel-ind-nf}~\ref{sel-ind-nf-a}. It leads to a short reproof of the \'{e}tale (or fppf) cohomological interpretation of Shafarevich--Tate groups and also forms the basis of our approach to fppf cohomological interpretation of Selmer groups.

\blem\lab{sel-loc-hens} Let $R$ be a discrete valuation ring, set $K\ce \Frac R$ and $K^h \ce \Frac R^h$, and let $\cG$ be a flat $R$-group algebraic space of finite presentation. If the horizontal arrows are injective in 
\[\xymatrix{
  H^1_\fppf(R, \cG) \ar@{^(->}[r] \ar[d] & H^1_\fppf(K, \cG_K) \ar[d] \\
  H^1_\fppf(R^h, \cG_{R^h}) \ar@{^(->}[r] & H^1_\fppf(K^h, \cG_{K^h}),
}\] 
then the square is Cartesian. If $\cG$ is a quasi-affine $R$-group scheme, then the same conclusion holds under analogous assumptions with $R^h$ and $K^h$ replaced by $\wh{R}$ and $\wh{K}$.
\elem

\bpf
We first treat the case of $R^h$ and $K^h$. We need to show that every $\cG_{K}$-torsor $\cT_{K}$ which, when base changed to $K^h$, extends to a $\cG_{R^h}$-torsor $\cT_{R^h}$, already extends to a $\cG$-torsor $\cT \ra \Spec R$. By \Cref{eq-cat}~\ref{eq-cat-b}, $\cT_{R^h}$ descends to a flat and of finite presentation $R$-algebraic space $\cT$, and various diagrams defining the $\cG$-action descend, too. To conclude that $\cT$ is a $\cG$-torsor, it remains to note that
\be\lab{tors-equiv} 
\cG \times_{R} \cT \ra \cT \times_{R} \cT, \quad (g, t) \mapsto (gt, t) 
\ee
is an isomorphism, as may be checked over $R^h$. 

In the similar proof for $\wh{R}$ and $\wh{K}$, to apply \Cref{eq-cat} one recalls that if $\cG$ is a quasi-affine scheme, then so are its torsors, see \cite{SP}*{\href{http://stacks.math.columbia.edu/tag/0247}{0247}}.
\epf

Let $S$ be a Dedekind scheme, let $K$ be its function field. As in \S\ref{glue}, to clarify analogies in \Cref{sel-cart}, we set $K_{S, s} \ce \Frac \cO_{S, s}$ for a nongeneric~$s \in S$.

\bcor\lab{sel-cart} 
Let $\cG$ be a flat closed $S$-subgroup scheme of an $S$-group scheme that is the N\'{e}ron model of its generic fiber. Then the square
\be\ba\lab{sel-cart-diag}\xymatrix{
  H^1_\fppf(S, \cG) \ar@{^(->}[r] \ar[d] & H^1_\fppf(K, \cG_K) \ar[d] \\
  \prod_s H^1_\fppf(\cO_{S,  s}, \cG_{\cO_{S, s}}) \ar@{^(->}[r] &\prod_s  H^1_\fppf(K_{S, s}, \cG_{K_{S, s}}),
}\ea\ee
is Cartesian (the products are indexed by the nongeneric $s\in S$), and similarly with $\cO_{S, s}$ and $K_{S, s}$ replaced by $\cO_{S, s}^h$ and $K_{S, s}^h$ (resp., $\wh{\cO}_{S, s}$ and $\wh{K}_{S, s}$ if $\cG \ra S$ is quasi-affine).
 \ecor

\bpf The indicated injectivity in \eqref{sel-cart-diag} results from \Cref{sub-ner-pf} and from the compatibility of the formation of the N\'{e}ron model with localization, henselization, and completion (see \cite{BLR90}*{\S1.2, Prop.~4 and \S7.2, Thm.~1 (ii)} for these compatibilities). By \Cref{sel-loc-hens}, the diagram
\[ \xymatrix{ \prod_s H^1_\fppf(\cO_{S, s}, \cG_{\cO_{S, s}}) \ar@{^(->}[r] \ar[d] & \prod_s H^1_\fppf(K_{S, s}, \cG_{K_{S, s}}) \ar[d] \\
  \prod_s H^1_\fppf(\cO_{ S, s}^h, \cG_{\cO_{S, s}^h}) \ar@{^(->}[r] & \prod_s H^1_\fppf(K_{S, s}^h, \cG_{K_{S, s}^h}) }
\]
is Cartesian, and likewise for $\wh{\cO}_{S, s}$ and $\wh{K}_{S, s}$. It remains to argue that \eqref{sel-cart-diag} is Cartesian.

We need to show that every $\cG_K$-torsor $\cT_K$ which extends to a $\cG_{\cO_{S, s}}$-torsor $\cT_{\cO_{S, s}}$ for every nongeneric $s\in S$, already extends to a $\cG$-torsor $\cT$ (these torsors are schemes, see the proof of \Cref{sub-ner-pf}). Since $\cT_K \ra \Spec K$ inherits finite presentation from $\cG_K$, for some open dense $U \subset S$ it spreads out to a $\cT_U \ra U$ which is faithfully flat, of finite presentation, has a $\cG_U$-action, and for which the analogue of \eqref{tors-equiv} over $U$ is bijective. Consequently, $\cT_U$ is a $\cG_U$-torsor.

To increase $U$ by extending $\cT_U$ over some $s \in S - U$, we spread out $\cT_{\cO_{S, s}}$ to a $\cG_W$-torsor $\cT_W$ over some open neighborhood $W \subset S$ of $s$. By \Cref{sub-ner-pf}, the torsors $\cT_U$ and $\cT_W$ are isomorphic over $U \cap W$, which permits us to glue them and to increase $U$. By iterating this process we arrive at the desired $U = S$.
\epf

\brems
\remi
The closed subgroup assumption on the flat $S$-group scheme $\cG$ is used only to deduce the indicated injectivity in \eqref{sel-cart-diag}. If one assumes instead that $\cG$ is commutative finite flat, then the injectivity follows from the valuative criterion of properness; consequently, \Cref{sel-cart} also holds for such $\cG$. For further extensions of \Cref{sel-cart}, see \cite{Ces14}*{7.2--7.4}.

\remi \lab{no-flat}
The flatness of $\cG$ is actually not needed for \Cref{sel-cart} to hold. To justify this, let $\wt{\cG}$ be the schematic image of $\cG_K$ in $\cG$, so that $\wt{\cG}$ is $S$-flat and a closed $S$-subgroup scheme of the same N\'{e}ron model. The formation of $\wt{\cG}$ commutes with flat base change, in particular, with base change to $\cO_{S, s}$, to $\cO_{S, s}^h$, or to $\wh{\cO}_{S, s}$. By \cite{Ces14}*{2.11} (or already by \cite{GMB13}*{Prop.~3.1} if $\cG$ is affine), the change of group maps
\[
\quad \quad \quad H^1_\fppf(S, \wt{\cG}) \ra H^1_\fppf(S, \cG) \quad \quad \text{and} \quad \quad H^1_\fppf(\cO_{S, s}, \wt{\cG}_{\cO_{S, s}}) \ra H^1_\fppf(\cO_{S, s}, \cG_{\cO_{S, s}})
\]
are bijective, and likewise with $\cO_{S, s}$ replaced by $\cO_{S, s}^h$ or by $\wh{\cO}_{S, s}$. This reduces the claim of \Cref{sel-cart} for $\cG$ to its claim for $\wt{\cG}$, which is $S$-flat.
\erems

We now use \Cref{sel-cart} to give an alternative proof of the results of \cite{Maz72}*{Appendix}.

\bprop\lab{mazur} Suppose that $S$ is a proper smooth curve over a finite field or that $S$ is the spectrum of the ring of integers of a number field. Let $A \ra \Spec K$ be an abelian variety, and let $\cA \ra S$ be its N\'{e}ron model. Letting the product run over the nongeneric $s\in S$, set
\[
\textstyle \Sha(\cA) \ce \Ker\p{H^1_\et(S, \cA) \ra \prod_{s } H^1_\et(\wh{\cO}_{S, s}, \cA_{\wh{\cO}_{S, s}})}.
\]
\benum
\item\lab{mazur-a} If $c_s$ denotes the local Tamagawa factor of $A$ at $s$ (see \S\ref{loc-tam} for the definition), then 
\[
\textstyle [H^1_\et(S, \cA) : \Sha(\cA)] \le \prod_{s } c_s.
\]

\item \lab{mazur-b} One has 
\[
\textstyle\Sha(\cA) = \Ker\p{H^1(K, A) \ra \prod_{s}  H^1(\wh{K}_{S, s}, A) }.
\]

\item One has
\[
\Sha(\cA) =   \im\p{H^1_\et(S, \cA^0) \ra H^1_\et(S, \cA)}.
\]

\item\lab{mazur-d} Let $\Sha(A)$ be the Shafarevich--Tate group of $A$. Then 
\[
\Sha(A) \subset \Sha(\cA)
\] 
and 
\[
\quad \quad [\Sha(\cA) : \Sha(A)] \le \prod_{\text{real } v} \#\pi_0(A(K_v)) \le 2^{\#\{\text{real } v\}\cdot \dim A}.
\]
In particular, $\Sha(A)$ is finite if and only if so is $H^1_\et(S, \cA)$.
\eenum
\eprop

\bpf \hfill
\benum 
\item By \Cref{nat-iso} (see \cite{Gro68}*{11.7 1$^\circ$)} for the identification between the \'{e}tale and the fppf cohomology groups), 
\[
\#H^1_\et(\wh{\cO}_{S, s}, \cA_{\wh{\cO}_{S, s}}) = c_s,
\]
so the claim results from the definition of $\Sha(\cA)$.

\item 
By \cite{BLR90}*{\S3.6, Cor.~10}, if an $A_{K_{S, s}^h}$-torsor has a $\wh{K}_{S, s}$-point, then it already has a $K_{S, s}^h$-point, i.e., the pullback map
\[
\quad \quad H^1(K_{S, s}^h, A) \ra H^1(\wh{K}_{S, s}, A)
\]
is injective, and hence, by \Cref{sub-ner-pf}, so is the pullback map
\[
\quad \quad H^1_{\et}(\cO_{S, s}^h, \cA_{\cO_{S, s}^h}) \ra H^1_{\et}(\wh{\cO}_{S, s}, \cA_{\wh{\cO}_{S, s}}).
\]
Therefore, it suffices to prove that
\[
\quad\quad \quad \textstyle \Ker\p{H^1_\et(S, \cA) \ra \prod_{s } H^1_\et(\cO^h_{S, s}, \cA_{\cO^h_{S, s}})} = \Ker\p{H^1(K, A) \ra \prod_{s}  H^1(K^h_{S, s}, A) }.
\]
This equality follows from the fact that the square
\[
\quad \xymatrix{
  H^1_\et(S, \cA) \ar@{^(->}[r] \ar[d] & H^1(K, A) \ar[d] \\
  \prod_s H^1_\et(\cO^h_{S, s}, \cA_{{\cO}^h_{S, s}}) \ar@{^(->}[r] &\prod_{s}  H^1(K^h_{S, s}, A)
}
\]
is Cartesian by \Cref{sel-cart}.

\item 
In the notation of \Cref{ses-gam-gam}, we have the exact sequence
\[
\quad \quad \tst 0 \ra \cA^0 \ra \cA \ra \bigoplus_s i_{s*} \Phi_s \ra 0.
\]
A segment of its associated long exact cohomology sequence reads
\[
\quad \quad \tst H^1_\et(S, \cA^0) \ra H^1_\et(S, \cA) \ra \bigoplus_s H^1_\et(k(s), \Phi_s),
\]
so it remains to recall that the pullback maps 
\[
\quad \quad H^1_\et(\wh{\cO}_{S, s}, \cA_{\wh{\cO}_{S, s}}) \ra H^1_\et(k(s), \Phi_s)
\]
are isomorphisms by \Cref{nat-iso}.

\item 
The inclusion follows from \ref{mazur-b}. So does the bound on the index because for real $v$ one has 
\[
\quad \quad H^1(K_v, A) \cong \pi_0(A(K_v)) \quad \quad \text{and} \quad \quad \#\pi_0(A(K_v)) \le 2^{\dim A},
\]
for instance, by \cite{GH81}*{1.1 (3) and 1.3}. The last claim also uses \ref{mazur-a}. \qedhere
\eenum
\epf

\section{Selmer groups as flat cohomology groups}\lab{selmer-flat}

The main goal of this section is the comparison of $\Sel_\phi A$ and $H^1_\fppf(S, \cA[\phi])$ in \Cref{sel-fppf-gl}.

\bpp[Selmer structures] \lab{sel-str} Let $K$ be a global field, and let $M$ be a finite discrete $\Gal(K^s/K)$-module. A \emph{Selmer structure} on $M$ is a choice of a subgroup of $H^1(K_v, M)$ for each place $v$ such that for all $v$ but finitely many, $H^1_\nr(K_v, M) \subset H^1(K_v, M)$ is chosen (compare with the definition \cite{MR07}*{1.2} in the number field case). The \emph{Selmer group} of a Selmer structure is the subgroup of $H^1(K, M)$ obtained by imposing the chosen local conditions, i.e., it consists of the cohomology classes whose restrictions to every $H^1(K_v, M)$ lie in the chosen subgroups. \epp

\bpp[The setup] If $K$ is a number field, we let $S \ce \Spec \cO_K$; if $K$ is a function field, we let $S$ be a connected proper smooth curve over a finite field with function field $K$. We let 
\[
A \xra{\phi} B \quad \quad \text{and} \quad \quad \cA \xra{\phi} \cB
\]
be a $K$-isogeny between abelian varieties and the induced $S$-homomorphism between their N\'{e}ron models. For a place $v\nmid \infty$, we get the induced map
\[
\phi_v\colon \Phi_{A, v} \ra \Phi_{B, v}
\]
between the groups of connected components of the special fibers of $\cA$ and $\cB$ at $v$. We let 
\[
c_{A, v} \ce \#\Phi_{A, v}(\bF_v) \quad \quad \text{and} \quad \quad  c_{B, v} \ce \#\Phi_{B, v}(\bF_v)
\]
be the local Tamagawa factors. \epp

\bpp[Two sets of subgroups (compare with \S\ref{3-subgp})]\lab{two-sel-str} 
The first set of subgroups is 
\[
\im\p{B(K_v) \xra{\kappa_{\phi, v}} H^1_\fppf(K_v, A[\phi])} \cong B(K_v)/\phi A(K_v) \subset H^1_\fppf(K_v, A[\phi]) \quad \quad \text{for all $v$.}
\]
Its Selmer group, defined as in \S\ref{sel-str}, is the \emph{$\phi$-Selmer group} 
\[
\Sel_\phi A \subset~H^1_\fppf(K, A[\phi]).
\]

The second set of subgroups is
\[\ba 
H^1_\fppf(\cO_v, \cA[\phi]) &\subset H^1_\fppf(K_v, A[\phi]),\quad\quad\text{ if }v\nmid \infty,\quad \text{ and} \\ 
H^1(K_v, A[\phi]) &\subset H^1(K_v, A[\phi]),\,\, \quad\quad\quad\text{if }v\mid \infty;
\ea\]
the indicated injectivity for $v\nmid \infty$ has been discussed in \S\ref{3-subgp} (even in the case when $\cA \xra{\phi} \cB$ fails to be flat!). By \Cref{sel-cart} and Remark \ref{no-flat} (together with \Cref{Aphi-aff}), its Selmer group is 
\[
H^1_\fppf(S, \cA[\phi]) \subset H^1_\fppf(K, A[\phi]).
\]

If $A[\phi]$ is \'{e}tale, then $\cA[\phi]$ is also \'{e}tale over a sufficiently small nonempty open subset of $S$, so, by \Cref{unr-comp}~\ref{unr-comp-e}, the above sets of subgroups are two Selmer structures on $A[\phi]$. 

In general, without assuming that $A[\phi]$ is \'{e}tale, the two sets of subgroups form two sets of Selmer conditions in the sense of \cite{Ces16}*{\S3.1}; in particular, by \cite{Ces16}*{3.2}, 
\[
H^1_\fppf(S, \cA[\phi]) \quad \quad \text{is always finite,}
\]
even in the case when $A[\phi]$ is not \'{e}tale and $\cA \xra{\phi} \cB$ is not flat. (The notion of \emph{Selmer conditions} generalizes the notion of a Selmer structure to the case when $M$ of \S\ref{sel-str} is an arbitrary commutative finite $K$-group scheme, i.e.,~not necessarily \'{e}tale.)
\epp

\bprop\lab{sel-fppf-gl} Suppose that $\cA \xra{\phi} \cB$ is flat (by \Cref{surj-id}, this assumption holds if, for example, $A$ has semiabelian reduction at all $v\nmid \infty$ with $\Char \bF_v \mid \deg \phi$).
\benum
\item\lab{sel-fppf-gl-c} 
If $\deg \phi$ is  prime to $\prod_{v\nmid \infty} c_{B, v}$, then 
\[
\Sel_\phi A \subset H^1_\fppf(S, \cA[\phi])
\]
inside $H^1_\fppf(K, A[\phi])$.

\item\lab{sel-fppf-gl-d} If $\deg \phi$ is 
 prime to $\prod_{v\nmid \infty} c_{A, v}$ and either $2 \nmid \deg \phi$ or $A(K_v)$ equipped with its archimedean topology is connected for all real $v$, then 
\[
H^1_\fppf(S, \cA[\phi]) \subset \Sel_\phi A
\] 
inside $H^1_\fppf(K, A[\phi])$.

\item\lab{sel-fppf-gl-e} If $\deg \phi$ is 
 prime to $\prod_{v\nmid \infty} c_{A, v}c_{B, v}$ and either $2\nmid \deg \phi$ or $A(K_v)$ equipped with its archimedean topology is connected for all real $v$, then
\[
H^1_\fppf(S, \cA[\phi]) = \Sel_\phi A
\]
inside $H^1_\fppf(K, A[\phi])$. 

\eenum
\eprop

\bpf 
By \S\ref{two-sel-str}, setting $H^1_\fppf(\cO_v, \cA[\phi]) \ce H^1(K_v, A[\phi])$ for $v\mid \infty$, we have injections
\be\lab{temp-lab}\ba 
\f{\Sel_\phi A}{H^1_\fppf(S, \cA[\phi])\cap \Sel_\phi A} &\hra \prod_{v\,\nmid\,\infty} \f{\im \kappa_{\phi, v}}{H^1_\fppf(\cO_v, \cA[\phi])\cap \im \kappa_{\phi, v}}, \\
\f{H^1_\fppf(S, \cA[\phi])}{H^1_\fppf(S, \cA[\phi])\cap \Sel_\phi A} &\hra \prod_v \f{H^1_\fppf(\cO_v, \cA[\phi])}{H^1_\fppf(\cO_v, \cA[\phi])\cap \im \kappa_{\phi, v}}.
\ea\ee
This together with \Cref{3-sub-comp}~\ref{3-sub-comp-d}, \ref{3-sub-comp-c}, and \ref{3-sub-comp-e} gives the claim because under the assumptions of \ref{sel-fppf-gl-d} and \ref{sel-fppf-gl-e} the factors of \eqref{temp-lab} for $v\mid \infty$ vanish: $H^1(K_v, A[\phi]) = 0$ unless $2\mid \deg \phi$ and $v$ is real, and also, by \cite{GH81}*{1.3}, $ H^1(K_v, A) \cong \pi_0(A(K_v))$.
 \epf

\brems 
\remi \lab{quant}
To compare $\Sel_{\phi} A$ and $H^1_\fppf(S, \cA[\phi])$ quantitatively, one may combine \eqref{temp-lab} with \Cref{3-sub-comp}~\ref{3-sub-comp-a}.

\remi 
As in \Cref{unr-comp}~\ref{unr-comp-f} and \ref{unr-comp-e}, the assumptions on $c_{A, v}$ and $c_{B, v}$ in \Cref{sel-fppf-gl}~\ref{sel-fppf-gl-c}, \ref{sel-fppf-gl-d}, and \ref{sel-fppf-gl-e} (and hence also in \Cref{sel-ind-nf}~\ref{sel-ind-nf-b}) can be weakened to, respectively,
\[\ba \#\Phi_{B, v}(\bF_v) &= \#(\phi_v(\Phi_{A, v}))(\bF_v) \text{ for all }v\nmid \infty,\\
 \#\Phi_{A, v}(\bF_v) &= \#(\phi_v(\Phi_{A, v}))(\bF_v) \text{ for all }v\nmid \infty, \text{ and}\\
 \#\Phi_{A, v}(\bF_v) &= \#(\phi_v(\Phi_{A, v}))(\bF_v) = \#\Phi_{B, v}(\bF_v) \text{ for all }v\nmid \infty.\ea\]
 
\remi \lab{rem-last}
In practice it is useful to not restrict \Cref{sel-fppf-gl} to the case when $A$ has semiabelian reduction at all $v\nmid \infty$ with $\Char \bF_v \mid \deg \phi$. For instance, suppose that $K$ is a number field, $A$ is an elliptic curve that has complex multiplication by an imaginary quadratic field $F \subset K$, and $\phi = \gA \in \End_K(A) \subset F \subset K$. Then 
\[
\quad \quad \cA_{\cO_K[\f{1}{\gA}]} \xra{\phi} \cA_{\cO_K[\f{1}{\gA}]}
\] 
is flat (even \'{e}tale) because it induces an automorphism of $\Lie \cA_{\cO_K[\f{1}{\gA}]}$, which is a line bundle on $\Spec \cO_K[\f{1}{\gA}]$. On the other hand, $\deg \phi$ need not be invertible on $\Spec \cO_K[\f{1}{\gA}]$. \Cref{sel-fppf-gl} applied to this example leads to a different proof of \cite{Rub99}*{6.4}, which facilitates the analysis of Selmer groups of elliptic curves with complex multiplication by relating them to class groups.
\erems

\appendix

\section{Torsors under a N\'{e}ron model}

\bpp[Dedekind schemes and N\'{e}ron models] \lab{ded-sch} 
A \emph{Dedekind scheme} $S$ is a connected Noetherian normal scheme of dimension $\le 1$. The connectedness is not necessary, but it simplifies the notation. We let $K$ denote the function field of $S$. An $S$-group scheme $\cX$ is a \emph{N\'{e}ron model} (of $\cX_K$) if it is separated, of finite type, smooth, and satisfies the \emph{N\'{e}ron property}: the restriction to the generic fiber map 
\[
\Hom_S(\cZ, \cX) \ra \Hom_K(\cZ_K, \cX_K)
\]
is bijective for every smooth $S$-scheme $\cZ$. 
\epp

\bprop\lab{tor-also-ner} Every torsor (for the fppf or the \'{e}tale topology) $\cT \ra S$ under a N\'{e}ron model $\cX \ra S$ is a scheme that is separated, smooth, and has the N\'{e}ron property. \eprop

\bpf 
Representability of $\cT$ by a scheme follows from \cite{Ray70}*{Thm.~XI 3.1 1)}. Its separatedness and smoothness are inherited from $\cX$ by descent. 

In checking the N\'{e}ron property, one can restrict to quasi-compact $\cZ$. Since $\cT$ is separated, $S$-morphisms $\cZ \xra{f} \cT$ are in bijection with closed subschemes 
\[
\mathfrak{Z} \subset \cZ \times_S \cT
\]
that are mapped isomorphically to $\cZ$ by the first projection ($\mathfrak{Z}$ is the graph of $f$), and similarly for $K$-morphisms $\cZ_K \ra \cT_K$. Such a $\mathfrak{Z}$ is determined by $\mathfrak{Z}_K$, being its schematic image in $\cZ \times_S \cT$ by \cite{EGAIV2}*{2.8.5}. Bijectivity of the assignment $\mathfrak{Z} \mapsto \mathfrak{Z}_K $ for any $\cZ$ as above is equivalent to the sought N\'{e}ron property of $\cT$. 

To check this bijectivity, it remains to show that the schematic image $\mathfrak{Z}\pr \subset \cZ \times_S \cT$ of any graph $\mathfrak{Z}_K \subset \cZ _K\times_K \cT_K$ is projected isomorphically to $\cZ$, as can be done \'{e}tale locally on $S$ (in the case of a Noetherian source, the formation of the schematic image commutes with flat base change by \cite{EGAIV3}*{11.10.3 (iv), 11.10.5 (ii)}). By \cite{EGAIV4}*{17.16.3 (ii)}, there is an \'{e}tale cover $S\pr \ra S$ trivializing the torsor $\cT$, so the claim follows from the N\'{e}ron property of $\cT_{S\pr} \cong \cX_{S\pr}$.
 \epf

\bcor \lab{inj-ner} For a N\'{e}ron model $\cX \ra S$, the pullback map
\be\lab{inj-ner-iso} H^1_\et(S, \cX) \xra{\iota} H^1_\et(K, \cX_K) {\cong} H^1(K, \cX_K)\ee
is injective. \ecor

\bpf 
By \Cref{tor-also-ner}, a torsor under $\cX$ is determined by its generic fiber. 
\epf

If $S$ is local, it is possible to determine the image of \eqref{inj-ner-iso}:

\bprop \lab{im-ner-unr} 
Suppose that $S = \Spec R$ for a discrete valuation ring $R$, and let $\cX \ra S$ be a N\'{e}ron model. The image of the injection $\iota$ from \eqref{inj-ner-iso} is the unramified cohomology subset
\[ I \ce \Ker(H^1(K, \cX_K) \ra H^1(K^{sh}, \cX_{K^{sh}}))\] 
where $K^{sh} \ce~\Frac R^{sh}$. In other words, an $\cX_K$-torsor $T$ extends to an $\cX$-torsor if and only if $T(K^{sh}) \neq \emptyset$. 
\eprop

\bpf 
Due to smoothness, every torsor $\cT$ under $\cX$ trivializes over an \'{e}tale cover $U \ra \Spec R$, and hence over $R^{sh}$, giving $\im \iota \subset I$. The inclusion $I \subset \im \iota$ is a special case of \cite{BLR90}*{\S6.5, Cor.~3}.
\epf

\Cref{inj-ner} can be strengthened as follows.

\bprop\lab{sub-ner-pf} For an $S$-flat closed $S$-subgroup scheme $\cG$ of a N\'{e}ron model $\cX \ra S$, the pullback map
\[ H^1_\fppf(S, \cG) \ra H^1_\fppf(K, \cG_K)\]
is injective. \eprop

\bpf In terms of descent data with respect to a trivializing $S\pr \ra S$ that is faithfully flat and locally of finite presentation, a $\cG$-torsor $\cT$ is described by the automorphism of the trivial right $\cG_{S\pr \times_S S\pr}$-torsor given by left translation by a $g \in \cG(S\pr \times_S S\pr)$. The image of $g$ in $\cX(S\pr \times_S S\pr)$ describes an $\cX$-torsor $\cT^{\cX}$, and the $\cG$-equivariant closed immersion $\cT \subset \cT^{\cX}$ of (a priori) algebraic spaces shows that $\cT$ is a scheme, since so is $\cT^{\cX}$ by \Cref{tor-also-ner}.

Let $\cT_1, \cT_2$ be $\cG$-torsors, and choose a common trivializing $S\pr \ra S$. It suffices to show that a $\cG_K$-torsor isomorphism $\gA_K\colon (\cT_1)_K \isomto (\cT_2)_K$ extends to a $\cG$-torsor isomorphism $\gA\colon \cT_1 \isomto \cT_2$. In terms of descent data, $\gA_K$ is described as left multiplication by a certain $h \in \cG(S\pr_K)$, whose image in $\cX(S\pr_K)$ extends $\gA_K$ to an $\cX_K$-torsor isomorphism $\beta_K\colon (\cT_1^{\cX})_K \isomto (\cT_2^\cX)_K$. By \Cref{tor-also-ner}, $\beta_K$ extends to an $\cX$-torsor isomorphism $\beta\colon \cT_1^\cX \isomto \cT_2^\cX$, which restricts to a desired $\gA$ due to schematic dominance considerations for $(\cT_i)_K \ra \cT_i$ (one uses \cite{EGAIV2}*{2.8.5} and \cite{EGAI}*{9.5.5}).
 \epf

\brem 
The above results continue to hold for N\'{e}ron lft models and without the flatness assumption in \Cref{sub-ner-pf}, see \cite{Ces14}*{2.19--2.21, 6.1} (an $S$-group scheme $\cX$ is a \emph{N\'{e}ron lft model} (of $\cX_K$) if it is separated, smooth, and satisfies the N\'{e}ron property recalled in \S\ref{ded-sch}; a N\'{e}ron lft model is not necessarily of finite type over $S$ but is always locally of finite type due to smoothness).
\erem

\section{Exact sequences involving N\'{e}ron models of abelian varieties}

In this appendix, we gather several standard facts about N\'{e}ron models of abelian varieties used in the main body of the paper.

\bpp[Open subgroups of N\'{e}ron models of abelian varieties]\lab{open} 
Let $S$ be a Dedekind scheme (defined in \S\ref{ded-sch}), and let $K$ be its function field. Let 
\[
A \ra \Spec K \quad \quad \text{and}\quad\quad \cA \ra S
\]
be an abelian variety and its N\'{e}ron model. For $s \in S$, let $\Phi_s \ce \cA_s/\cA_s^0$ be the \'{e}tale $k(s)$-group scheme of connected components of $\cA_s$.  For each nongeneric $s \in S$, choose a $k(s)$-subgroup $\Gamma_s \subset \Phi_s$. Then for all $s$ but finitely many, $\Gamma_s = \Phi_s$, and we define the open subgroup 
\[
\cA^\Gamma \subset \cA
\]
by removing for every $s$ the connected components of $\cA_s$ not in $\Gamma_s$. Letting $i_s\colon \Spec k(s) \ra S$ denote the inclusion of the nongeneric point $s$, we have the homomorphism 
\[
\textstyle\cA^{\Gamma} \ra \bigoplus_s i_{s*}\Gamma_s.
\]
If $\Gamma_s = 0$ for every $s$, then the resulting $\cA^0$ is the fiberwise identity component of $\cA$.
\epp

\bprop\lab{ses-gam-gam} For all choices $\wt{\Gamma}_s \subset \Gamma_s \subset \Phi_s$, the sequence 
\[ 
\textstyle 0 \ra \cA^{\wt{\Gamma}} \ra \cA^\Gamma \xra{a} \bigoplus_s i_{s*}(\Gamma_s/\wt{\Gamma}_s) \ra 0 
\]
is exact in $S_{\et}$, $S_\Et$, and $S_\fppf$. \eprop 

\bpf 
Left exactness is clear, whereas to check the remaining surjectivity of $a$ in $S_\Et$ on stalks, it suffices to consider strictly local $(\cO, \fm)$ centered at a nongeneric $s \in S$ with $\wt{\Gamma}_s \neq \Gamma_s$. Let $\fa \subset \fm$ be the ideal generated by the image of $\fm_{S, s}$. In the commutative diagram
\[
\xymatrix{
\cA^{\Gamma}(\cO) \ar[r]^-{a(\cO)} \ar@{->>}[d]^-{b} & (\Gamma_s/\wt{\Gamma}_s)(\cO/\fa) \ar[d]^-{d}_-{\wr} \\
\cA^{\Gamma}(\cO/\fm) \ar@{->>}[r]^-{c}  & (\Gamma_s/\wt{\Gamma}_s)(\cO/\fm),
}
\]
the surjectivity of $b$ follows from Hensel-lifting for the smooth $\cA^\Gamma_\cO \ra \Spec \cO$ (see \cite{EGAIV4}*{18.5.17}), the surjectivity of $c$ follows from the invariance of the component group of the smooth $\cA^{\Gamma}_{k(s)^s} \ra \Spec k(s)^s$ upon passage to a separably closed overfield, whereas the bijectivity of $d$ is immediate from $(\Gamma_s/\wt{\Gamma}_s)_{\cO/\fa}$ being finite \'{e}tale over the Henselian local $(\cO/\fa, \fm/\fa)$. The desired surjectivity of $a(\cO)$ follows.
 \epf

Let $A \xra{\phi} B$ be a $K$-isogeny of abelian varieties, and let $\cA \xra{\phi} \cB$ be the homomorphism induced on N\'{e}ron models over $S$.

\bprop \lab{Aphi-aff}
The kernel $\cA[\phi] \ra S$ is affine; every fppf torsor under $\cA[\phi]$ is representable.
\eprop

\bpf
Affineness of $\cA[\phi]$ is a special case of \cite{Ana73}*{2.3.2}. Effectivity of fppf descent for affine schemes gives the torsor claim.
\epf

\blem\lab{surj-id} The following are equivalent:
\benum
\item\lab{surj-id-c} $\cA \xra{\phi} \cB$ is quasi-finite,
\item\lab{surj-id-b} $\cA^0 \xra{\phi} \cB^0$ is surjective (as a morphism of schemes),
\item\lab{surj-id-a} $\cA \xra{\phi} \cB$ is flat,
\eenum
and are implied by 
\benum \addtocounter{enumi}{3}
\item\lab{surj-id-d} $A$ has semiabelian reduction at all the nongeneric $s \in S$ with $\Char k(s) \mid \deg \phi$.
\eenum
 \elem

\bpf 
Due to the fibral criterion of flatness \cite{EGAIV3}*{11.3.11} for \ref{surj-id-a}, the conditions \ref{surj-id-c}--\ref{surj-id-a} can be checked fiberwise on $S$. We will show that they are equivalent for the fiber over an $s \in S$.

Since $\cA$ and $\cB$ are faithfully flat and locally of finite type over $S$, \cite{BLR90}*{\S2.4, Prop.~4} supplies the equalities 
\[
\dim \cA_s = \dim A \quad \quad \text{and} \quad \quad \dim \cB_s = \dim B,
\] 
and hence also $\dim \cA_s = \dim \cB_s$. Moreover, by \cite{SGA3Inew}*{VI$_{\text{A}}$, 6.7}, every homomorphism between algebraic groups over a field factors through a flat surjection onto its closed image, so $\phi_s$ is surjective on identity components if and only if it is quasi-finite, i.e, \ref{surj-id-c}$\Leftrightarrow$\ref{surj-id-b}. Furthermore, if $\phi_s( \cA^0_s) = \cB_s^0$, then $\phi_s$ is flat on identity components, i.e.,~\ref{surj-id-b}$\implies$\ref{surj-id-a}. Conversely, if $\phi_s$ is flat, then, in addition to being closed, $\phi_s( \cA^0_s)$ is also open, and hence equals $\cB^0_s$, i.e.,~\ref{surj-id-a}$\implies$\ref{surj-id-b}.

For the last claim, the consideration of the isogeny $\psi\colon B \ra A$ with the kernel $\phi(A[\deg \phi])$ reduces to the case when $\phi$ is multiplication by an integer $n$. For such $\phi$, the surjectivity of $\phi_s$ on the identity components is clear if the reduction at $s$ is semiabelian and follows by inspection of Lie algebras if $\Char k(s) \nmid n$. 
\epf

\bcor\lab{qf-flat} 
Suppose that $\cA \xra{\phi} \cB$ is flat (e.g., that $A$ has semiabelian reduction at every nongeneric $s \in S$ with $\Char k(s) \mid \deg \phi$). Then $\cA[\phi] \ra S$ is quasi-finite, flat, and affine; it is also finite if $A$ has good reduction at every nongeneric point of $S$.
\ecor

\bpf By \Cref{surj-id}, $\cA \xra{\phi} \cB$ is quasi-finite and flat; in the good reduction case, it is finite due to its properness, see \cite{EGAIV3}*{8.11.1}. Affineness results from \Cref{Aphi-aff}.
 \epf

\bcor\lab{Aphi-ner} If $\Char k(s) \nmid \deg \phi$ for all $s \in S$, then $\cA[\phi]$ is the N\'{e}ron model of $A[\phi]$. \ecor

\bpf Due to \Cref{qf-flat} and the degree hypothesis, the quasi-finite flat $\cA[\phi] \ra S$ is \'{e}tale. On the other hand, by \cite{BLR90}*{\S7.1, Cor.~6}, the N\'{e}ron model of $A[\phi]$ may be obtained as the group smoothening of the schematic image of $A[\phi]$ in $\cA$. By \cite{EGAIV2}*{2.8.5}, this schematic image is $\cA[\phi]$, so, since $\cA[\phi] \ra S$ is \'{e}tale, no smoothening is needed. \epf

A choice of $k(s)$-subgroups $\Gamma_s \subset \Phi_s$ gives rise to their images $\phi_s(\Gamma_s)$. These images, in turn, give rise to the open subgroup $\cB^{\phi(\Gamma)} \subset \cB$ as in \S\ref{open}.

\bcor\lab{ses-mult} Suppose that $\cA \xra{\phi} \cB$ is flat (e.g., that $A$ has semiabelian reduction at all the nongeneric $s\in S$ with $\Char k(s) \mid \deg\phi$). Then for every choice of $k(s)$-subgroups $\Gamma_s \subset \Phi_s$, the sequence 
\[
\textstyle 0 \ra \cA^\Gamma[\phi] \ra \cA^\Gamma \xra{\phi} \cB^{\phi(\Gamma)} \ra 0
\]
is exact in $S_\fppf$. \ecor

\bpf The $S$-morphism $\cA^\Gamma \xra{\phi} \cB^{\phi(\Gamma)}$ is faithfully flat and locally of finite presentation by \Cref{surj-id}, whereas the exactness at the other terms is immediate from the definitions.
 \epf

\end{document}